\documentclass[twocolumn]{autartmod}

\usepackage[utf8]{inputenc}
\usepackage[T1]{fontenc}
\usepackage{mathpazo}
\usepackage{graphicx}
\usepackage{placeins}
\usepackage{amsmath,amssymb,bm}
\usepackage{xcolor}
\usepackage{silence}
\usepackage{subcaption}

\makeatletter
\renewcommand*\l@subsection{\@dottedtocline{2}{2.5em}{2.3em}} 
\makeatother

\allowdisplaybreaks 
\usepackage{ragged2e} 

\usepackage{tikz}
\usetikzlibrary{patterns} 
\usetikzlibrary{calc,positioning,shapes.geometric} 
\usepackage{pgfplots}
\pgfplotsset{compat=newest}
\usetikzlibrary{arrows.meta,pgfplots.fillbetween}
\newcommand{\coloredblock}[3][0.2cm]{\tikz[baseline=-0.75ex]{\node[fill=#2, draw=#3, minimum size=#1, inner sep=0pt] (n) {};}}

\usepackage{algorithm}
\usepackage{algpseudocode}

\colorlet{fadedred}{red!30!gray!80!white}
\colorlet{fadedblue}{blue!70!gray!80!white}
\colorlet{fadedgreen}{green!60!blue!80!gray}

\definecolor{beige}{RGB}{245,245,220}

\usetikzlibrary{shadows}
\usepackage[many]{tcolorbox}
\newtcolorbox{eqbox}{
  colback=beige!30,     
  colframe=black!30,    
  boxrule=0.7pt,        
  arc=5pt,              
  boxsep=0pt,           
  left=0pt, right=5pt,  
  top=-6pt, bottom=-10pt
}
\newtcolorbox{eqboxb}{
  colback=beige!30,     
  colframe=black!30,    
  boxrule=0.7pt,        
  arc=5pt,              
  boxsep=0pt,           
  left=0pt, right=5pt,  
  top=-6pt, bottom=5pt
}

\usepackage{hyperref}
\hypersetup{
    colorlinks=true,
    linkcolor=fadedblue, 
    citecolor=fadedgreen,
    urlcolor=fadedgreen
}

\newcommand{\behcet}{Beh\c{c}et}
\newcommand{\acikmese}{A\c{c}\i kme\c{s}e}

\renewenvironment{pf}{\par\textcolor{gray}{\textit{\textbf{Proof:}~}}}{\hfill\coloredblock{black!50}{black!50}\par}

\let\ae\undefine
\newcommand{\ae}{\mathrm{a.e.}}
\newcommand{\Neq}{n_{h}}
\newcommand{\Nineq}{n_{g}}
\newcommand{\tinit}{t_{\mathrm{i}}}
\newcommand{\tfinal}{t_{\mathrm{f}}}
\newcommand{\tspan}{[\tinit,\tfinal]}

\newcommand{\rinit}{r_{\mathrm{i}}}
\newcommand{\rfinal}{r_{\mathrm{f}}}
\newcommand{\derv}[1]{\overset{{\scriptscriptstyle\circ}}{#1}}
\renewcommand{\dot}[1]{\overset{\text{{\large{.}}}}{#1}}
\newcommand{\penaltyfun}{\Lambda}
\newcommand{\bR}{\mathbb{R}}
\newcommand{\bRnx}{\bR^{n_x}}
\newcommand{\bRnu}{\bR^{n_u}}
\newcommand{\bRp}{\bR_{\scriptscriptstyle +}}
\newcommand{\mc}[1]{\mathcal{#1}}
\newcommand{\mb}[1]{\mathbb{#1}}
\newcommand{\indic}{\tilde{I}}
\newcommand{\ncone}{\mc{N}}
\newcommand{\dbyd}[2]{\frac{\mathrm{d}#1}{\mathrm{d}#2}}

\newcommand{\inlpbyp}[2]{\nabla_{#2}#1}
\newcommand{\normplus}[1]{|#1|_{\scriptscriptstyle +}}
\newcommand{\eye}[1]{I_{#1}}
\newcommand{\ones}[1]{{1}_{#1}}
\newcommand{\zeros}[1]{{0}_{#1}}
\newcommand{\selector}[1]{{{}^{#1}{\!E\;\!}}}
\newcommand{\ntildx}{n_{\tilde{x}}}
\newcommand{\nobs}{n}
\newcommand{\scvxgen}{\textsc{sc{\scalebox{0.675}{vx}}gen}}

\usepackage{comment}
\usepackage{todonotes}
\edef\endfrontmatter{%
  \unexpanded\expandafter{\endfrontmatter}
  \noexpand\endNoHyper 
}

\begin{document}

\begin{frontmatter}
\runtitle{SCP Traj. Opt. CTCS}

\title{%
Successive Convexification for Trajectory Optimization with\\Continuous-Time Constraint Satisfaction\!%
\thanksref{footnoteinfo}}

\thanks[footnoteinfo]{\!\!Corresponding author: Purnanand Elango.}



\author{Purnanand Elango}$^{*}$\ead{pelango@uw.edu},~~%
\author{Dayou Luo}$^{\dagger}$\ead{dayoul@uw.edu},~~%
\author{Abhinav G.\ Kamath}$^{*}$\ead{agkamath@uw.edu},~~%
\author{Samet Uzun}$^{*}$\ead{samet@uw.edu},
\author{Taewan Kim}$^{*}$\ead{twankim@uw.edu},~~and~~%
\author{\behcet~\acikmese}$^{*}$\ead{behcet@uw.edu}%

\address{$^{*}$William E.\ Boeing Department of Aeronautics \& Astronautics}
\address{$^{\dagger}$Department of Applied Mathematics}
\address{\vphantom{$^{\dagger}$}University of Washington, Seattle, WA, 98195}


\begin{abstract} 
We present successive convexification, a real-time-capable solution method for nonconvex trajectory optimization, with continuous-time constraint satisfaction and guaranteed convergence, that only requires first-order information. The proposed framework combines several key methods to solve a large class of nonlinear optimal control problems: (i) exterior penalty-based reformulation of the path constraints; (ii) generalized time-dilation; (iii) multiple-shooting discretization; (iv) $\ell_1$ exact penalization of the nonconvex constraints; and v) the prox-linear method, a sequential convex programming (SCP) algorithm for convex-composite minimization. The reformulation of the path constraints enables continuous-time constraint satisfaction even on sparse discretization grids and obviates the need for mesh refinement heuristics. Through the prox-linear method, we guarantee convergence of the solution method to stationary points of the penalized problem and guarantee that the converged solutions that are feasible with respect to the discretized and control-parameterized optimal control problem are also Karush-Kuhn-Tucker (KKT) points. Furthermore, we highlight the specialization of this property to global minimizers of convex optimal control problems, wherein the reformulated path constraints cannot be represented by canonical cones, i.e., in the form required by existing convex optimization solvers. In addition to theoretical analysis, we demonstrate the effectiveness and real-time capability of the proposed framework with numerical examples based on popular optimal control applications: dynamic obstacle avoidance and rocket landing.
\end{abstract}

\end{frontmatter}


\section{Introduction}
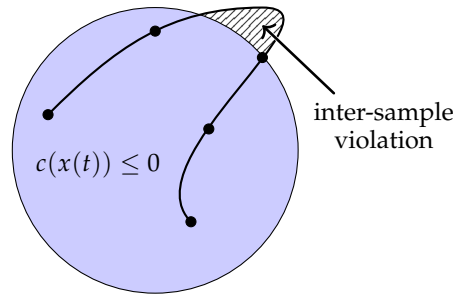
\begin{figure}[!htpb]
    \centering


\begin{tikzpicture}[scale=0.95, transform shape, line cap=round]

\begin{scope}
    \fill[pattern=north east lines, pattern color=darkgray] (-1.5,0.5) to[out=45,in=180] (1.5,2) to[out=0,in=135] (0.5,-1) -- cycle;
\end{scope}

\filldraw[fill=blue!20,draw=black] (0,0) ellipse (2cm);

\node at (-0.8,-0.2) {$c(x(t)) \le 0$};

\draw[thick] (-1.5,0.5) to[out=45,in=180] (1.5,2) to[out=0,in=135] (0.5,-1);

\filldraw (0.5,-1) circle (2pt);
\filldraw (-1.5,0.5) circle (2pt);
\filldraw (1.5,1.3) circle (2pt);
\filldraw (0,1.67) circle (2pt);
\filldraw (0.75,0.3) circle (2pt);

\node[rectangle, draw=white, fill=white, text width=2cm, align=center] at (3.2,0.35) {inter-sample \\ violation};
\draw[thick, <-, draw=white, line width=1.5pt] (1.5, 1.75) -- (2.5,0.8);
\draw[thick, <-, draw=black, line width=1pt] (1.5, 1.75) -- (2.5,0.8);
  
\end{tikzpicture}

    \caption{Direct methods for trajectory optimization impose path constraints (such as $c(x(t))\le 0$) at finitely-many time nodes (black dots), invariably leading to inter-sample violation (shaded region) in the state trajectory $x$ obtained through simulation of the dynamical system with the control input solution. The proposed framework mitigates this phenomenon.}
    \label{fig:motivate}
\end{figure}
Trajectory optimization forms an important part of modern guidance, navigation, and control (GNC) systems, wherein, it is used to generate reference trajectories for onboard use and also in offline design and analysis. State-of-the-art trajectory optimization methods do not yet check all the boxes with regard to desirable features: continuous-time feasibility, real-time performance, convergence guarantees, and numerical robustness \cite{malyuta2021advances}. Among these, continuous-time feasibility, which refers to the feasibility of the state trajectory and the control input with respect to the system dynamics \textit{and} path constraints, is especially challenging due to the infinite-dimensional nature of optimal control problems. However, continuous-time feasibility is essential for meeting safety and performance requirements, which is a prerequisite for the deployment of autonomous systems such as space exploration vehicles, reusable rockets, dexterous robotic manipulators, etc \cite{malyuta2022convex}. Aerospace applications, in particular, typically prioritize feasibility and real-time performance over optimality \cite{tsiotras2017toward}.

Trajectory optimization methods can be broadly classified into \textit{indirect} and \textit{direct methods} \cite[Sec. 4.3]{betts2010practical}. Invoking the maximum principle to solve optimal control problems with general nonlinear dynamics subject to path constraints is challenging. Such solution methods, i.e., indirect methods, usually take an \textit{optimize-then-discretize} approach, and can only handle a restricted class of problems. The solutions they generate, however, naturally ensure continuous-time feasibility \cite{betts1998survey,rao2009survey}. Methods that take the \textit{discretize-then-optimize} approach, i.e., direct methods, on the other hand, are capable of handling general optimal control problems and are more reliable numerically (they are less sensitive to initialization). Since discretization is a crucial step in direct methods, the resulting solutions invariably suffer from so-called inter-sample constraint violations \cite{dueri2017trajectory} (see Figure \ref{fig:motivate}).

A majority of the existing direct methods, such as MISER \cite{jennings1991miser3}, DIRCOL \cite{vonStryk1993dircol}, 
PSOPT \cite{becerra2010psopt}, PROPT \cite{rutquist2010propt}, and GPOPS-II \cite{patterson2014gpops2}, transcribe optimal control problems to nonlinear programs (NLP) and call standard NLP solvers, such as SNOPT \cite{gill2006snopt}, IPOPT \cite{wachter2006ipopt}, and \textsc{Knitro} \cite{byrd2006knitro}, which require second-order information, lack convergence guarantees, and are unsuitable for real-time, embedded applications.  In contrast, software such as \texttt{acados} \cite{verschueren2022acados} and ACADO \cite{houska2011acado} interface with custom quadratic programming (QP), sequential quadratic programming (SQP), and NLP solvers that exploit the structure of the underlying problem. However, they neither guarantee convergence nor enforce continuous-time constraint satisfaction.

Several direct methods either---(i) assume a discrete-time formulation with path constraints imposed at finitely-many time nodes, i.e., they disregard the conversion of continuous-time problem descriptions to discrete-time ones (e.g., TrajOpt \cite{schulman2014motion}, PANOC \cite{stella2017simple,sathya2018embedded}, CALIPSO \cite{howell2023calipso}); (ii) provide convergence guarantees by restricting the class of path constraints \cite{liu2018convex,roulet2019iterative,lee2022convexifying,zhang2021optimization,liu2014solving,marcucci2022motion}; or (iii) ensure continuous-time constraint satisfaction for a  restricted class of dynamical systems \cite{acikmese2006convex,acikmese2008enhancements,richards2015intersample,dueri2017trajectory}. 

\textit{Sequential convex programming} (SCP) algorithms for trajectory optimization have received attention in the recent years \cite{messerer2021survey,malyuta2022convex} as a competitive alternative to sequential quadratic programming (SQP) and IPM-based NLP algorithms, especially since SCP is a multiplier-free method and does not require second-order information. In trajectory optimization, computing the second-order sensitivities of nonlinear dynamics (which form a part of the Hessian of the Lagrangian within SQP and IPM) is expensive \cite{diehl2022book}. Recent work has explored convergence guarantees for SCP, in the context of both trajectory optimization \cite{liu2014solving,foust2020optimal,bonalli2019gusto,bonalli2022analysis,xie2023descent,xie2023higher} and general nonconvex optimization \cite{dinh2010local,messerer2020determining,lewis2016proximal,cartis2011evaluation,drusvyatskiy2018error}. SCP-based trajectory optimization methods, without theoretical guarantees, have been applied to a wide range of robotics and aerospace applications, with domain-specific heuristics that ensure effective practical performance \cite{schulman2014motion,wang2017constrained,szmuk2020successive,sagliano2018pseudospectral}. Widespread adoption of SCP for performance- and safety-critical applications, however, will necessitate the development of a general framework with rigorous certification of its capabilities.

We propose \textit{successive convexification}, an SCP-based, real-time-capable solution method for nonconvex trajectory optimization, with continuous-time feasibility and guaranteed convergence. The proposed framework combines several key methods to solve a large class of nonlinear optimal control problems: (i) exterior penalty-based reformulation of the path constraints; (ii) generalized time-dilation; (iii) multiple-shooting discretization; (iv) $\ell_1$ exact penalization of the nonconvex constraints; and (v) the prox-linear method, an SCP algorithm for convex-composite minimization. 

The reformulation of path constraints involves integrating the continuous-time constraint violation using a smooth exterior penalty, which is transformed into an auxiliary dynamical system with boundary conditions. The reformulation combined with multiple-shooting discretization \cite{bock1984multiple} enables continuous-time feasibility even on sparse discretization grids, and obviates the need for mesh refinement heuristics. While such constraint reformulations have appeared in the optimal control literature since the 1960s \cite{mcgill1965optimum,sargent1977development,teo1987simple} with specific choices of penalty functions \cite{loxton2009optimal,lin2014exact}, the full extent of its capabilities for enabling SCP-based convergence-guaranteed, continuous-time-feasible trajectory optimization have not been explored, to the best of the authors' knowledge. We address the consequences of the reformulation on a constraint qualification that is important in numerical optimization, and provide a rigorous quantification of the extent of continuous-time constraint satisfaction. 

Through the prox-linear method \cite{drusvyatskiy2018error,drusvyatskiy2019efficiency}, we guarantee convergence of the solution method to stationary points of the $\ell_1$-penalized problem and guarantee that the converged solutions that are feasible with respect to the discretized and control-parameterized optimal control problem are also Karush-Kuhn-Tucker (KKT) points. Furthermore, we highlight the specialization of this property to global minimizers of convex optimal control problems, wherein the reformulated path constraints cannot be represented by canonical cones, i.e., in the form required by existing convex optimization solvers. 

In addition to theoretical analysis, we demonstrate the effectiveness of the proposed framework with numerical examples based on popular optimal control applications: two nonconvex problems---dynamic obstacle avoidance and 6-DoF rocket landing, and one convex problem---3-DoF rocket landing using lossless convexification. We also demonstrate the real-time capability of the proposed framework on the nonconvex examples considered, by executing C code generated using {\scvxgen}, an in-house-developed general-purpose real-time trajectory optimization software with customized code-generation support. The C codebase generated by {\scvxgen} uses the proposed framework to solve the optimal control problem at hand. 

Simplified implementations of the proposed framework were recently demonstrated for specific applications, ranging from GPU-accelerated trajectory optimization for six-degree-of-freedom (6-DoF) powered-descent guidance \cite{chari2024fast} and nonlinear model predictive control (NMPC) for obstacle avoidance \cite{uzun2024nmpc}, to trajectory optimization for 6-DoF aircraft approach and landing \cite{kim2024approach}.

\section{Problem Formulation}
This section describes the transformation of a path-constrained, free-final-time optimal control problem to a fixed-final-time optimal control problem through generalized time-dilation and constraint reformulation. The constraint reformulation involves the conversion of path constraints into a two-point boundary value problem for an auxiliary dynamical system.
%
\subsection{Notation}
We adopt the following notation in the remainder of the discussion. %
The set of real numbers is denoted by $\bR$, the set of nonnegative real numbers by $\bRp$, the set of real $n\times m$ matrices by $\bR^{n\times m}$, and the set of real $n\times 1$ vectors by $\bR^n$. %
The concatenation of vectors $v\in\bR^n$ and $u\in\bR^m$ is denoted  by $(v,u)\in\bR^{n+m}$, the concatenation of matrices $A\in\bR^{l\times m}$ and $B\in\bR^{l\times n}$ by $[A~B] \in \bR^{l\times(m+n)}$, the Cartesian product of sets $C$ and $D$ by $C\times D$, and the Kronecker product by $\otimes$. %
The vector of ones in $\bR^n$ is denoted by $\ones{n}$, the identity matrix in $\bR^{n\times n}$ by $\eye{n}$, and the matrix of zeros in $\bR^{n\times m}$ by $\zeros{n\times m}$. Whenever the subscript is omitted, the size is inferred from context. %
For any scalar $v$, we define $\normplus{v} = \max\{0,v\}$. The operations $\normplus{\square}$, $|\square|$, $\square^2$ (and their compositions) apply elementwise for a vector. %
The Euclidean norm of a vector $v$ is denoted by $\|v\|$. %
The indicator function of a convex set $D$ is denoted by $\indic_D$ (see \cite[E.g. 3.1]{boyd2004convex}), and its normal cone at $x$ by $\ncone_D(x)$ \cite[Def. 5.2.3]{hiriart-urruty1993convex}. %
The subdifferential of a function $f$, evaluated at $x$, is a set denoted by $\partial{f}(x)$, and its members are called subgradients. %
The gradient of a differentiable function $g:\bR^{n}\to\bR^l$ with respect to $z\in\bR^m$, evaluated at $x\in\bR^n$, is denoted by $\nabla_z g(x)\in\bR^{n\times m}$, where the elements of $z$ can include a subset or superset of the arguments of $g$ (irrespective of the order), and the subscript $z$ is omitted if it coincides with the list of arguments of $g$. The (sub)gradient of a scalar-valued function at a point is defined to be a row vector (for e.g., from $\bR^{1\times m}$). The (sub)gradient of a vector-valued function is a matrix consisting of the (sub)gradients of its scalar-valued elements along the rows. In particular, the partial derivatives of function $(u,v)\mapsto h(u,v)$, evaluated at $\bar{u},\bar{v}$, are denoted by $\inlpbyp{h}{u}(\bar{u},\bar{v})$, $\inlpbyp{h}{v}(\bar{u},\bar{v})$, respectively (with the argument inferred from context whenever omitted). Furthermore, the notation $\dbyd{\square(\varsigma)}{\varsigma}$ represents the derivative with respect to a scalar variable.
%
\subsection{Optimal Control Problem}\label{subsec:ocp}
We consider a class of \textit{free-final-time} optimal control problems for nonlinear dynamical systems with nonconvex constraints on the state and input, given by
\begin{eqbox}
\begin{subequations}
\begin{align}
\underset{x,\,u,\,\tfinal}{\operatorname{minimize}}~\,&~L(\tfinal,x(\tfinal)) \\
\operatorname{subject~to}\,&~\dot{x}(t) = f(t,x(t),u(t))~\,\ae~t\in\tspan\label{ocp:dyn}\\
 &~g(t,x(t),u(t)) \le \zeros{n_g}\:\,\,\,\,\ae~t\in\tspan 
 \label{ocp:ineq-cnstr}\\
 &~h(t,x(t),u(t)) = \zeros{n_h}\:\,\,\,\,\ae~t\in\tspan \label{ocp:eq-cnstr}\\
 &~P(\tinit,x(\tinit),\tfinal,x(\tfinal)) \le \zeros{n_P} \label{ocp:bc-ineq}\\ 
 &~Q(\tinit,x(\tinit),\tfinal,x(\tfinal)) = \zeros{n_Q}\label{ocp:bc-eq} 
\end{align}\label{ocp}%
\end{subequations}
\end{eqbox}
where the derivative with respect to time $t\in\tspan$ in \eqref{ocp:dyn} is denoted by $\dot{\square} = \dbyd{\square}{t}$. We assume that the terminal state cost function $L:\bRp\times\bRnx\to\bR$, the dynamics function $f:\bRp\times\bRnx\times\bRnu\to\bRnx$, the path constraint functions $g:\bRp\times\bRnx\times\bRnu\to\bR^{\Nineq}$, $h:\bRp\times\bRnx\times\bRnu\to\bR^{\Neq}$, and the boundary condition constraint functions $P:\bRp\times\bRnx\times \bRp\times\bRnx\to\bR^{n_P},~Q:\bRp\times\bRnx\times \bRp\times\bRnx\to\bR^{n_Q}$, are continuously differentiable. The inequality and equality constraints in \eqref{ocp:ineq-cnstr}-\eqref{ocp:bc-eq} are interpreted elementwise. The initial time $\tinit\in\bRp$ is fixed, while the final time $\tfinal\in\bRp$ is a free (decision) variable. For simplicity, we omit parameters \cite[Eq. 1]{malyuta2022convex} of the system dynamics and constraints in \eqref{ocp}. With minimal modifications, the subsequent development can handle parameters as decision variables as well. Note that the above problem can be a \textit{nonconvex trajectory optimization} problem due to: i) nonlinear dynamics, i.e., $f$ is a nonlinear function; ii) final time being free; iii) nonlinear functions $h$ and $Q$; iv) nonconvex functions $g$ and $P$.  

We assume that the \textit{control input} $u:\tspan\to\bRnu$ is a \textit{piecewise continuous function}.  Since $f$ is continuously differentiable, this  ensures the existence and uniqueness of an absolutely continuous function $x : \tspan \to \bRnx$,  the \textit{state trajectory}, which satisfies \eqref{ocp:dyn} almost everywhere and the following integral equation
$$
    x(t) = x_{\mathrm{i}} + \int_{\tinit}^{t}f(\gamma,x(\gamma),u(\gamma))\mathrm{d}\gamma
$$
where $x(\tinit) = x_{\mathrm{i}}$ is the initial state. Note that existence and uniqueness of the state trajectory can be ensured with weaker assumptions; we refer the reader to \cite[Sec. 3.3.1]{liberzon2011calculus} and \cite[Sec. II.3]{berkovitz1974optimal} for detailed discussions. We also  require the following boundedness assumption based on Gronwall's Lemma \cite[Chap. 4, Prop. 1.4]{clarke1998nonsmooth}.
\begin{assum}\label{asm:x-bnd}
For any compact $\mb{U} \subset \bRnu$ and $\tfinal > \tinit \ge 0$, there exist positive  $\vartheta$ and $\varrho$ such that $\|f(t, x, u)\| \leq \vartheta \|x\| + \varrho$
for all $(t, x, u) \in [\tinit, \tfinal] \times \bRnx \times \mb{U}$.
\end{assum}

Note that Assumption \ref{asm:x-bnd} may not always hold. Consider $f(t,x,u) = x^2$. For any positive $\vartheta$ and $\varrho$, when $x$ becomes sufficiently large, Assumption \ref{asm:x-bnd} is invalid. Nonetheless, for physical systems, one can define a compact set $\mathbb{X}\subset\bRnx$, containing all physically meaningful states. Any state not contained in $\mathbb{X}$ can be regarded as infeasible. Consequently, we create a modified dynamics function so that it satisfies Assumption \ref{asm:x-bnd} and is consistent with the original dynamics function on $\mathbb{X}$. To be precise, for any open set $\mathbb{W}$ containing $\mathbb{X}$, there exists a smooth function $\omega : \bRnx \rightarrow \bR$ such that $\omega(x) = 1$ for all $x\in\mathbb{X}$ and the support of $\omega$ is compact and contained in $\mathbb{W}$, as shown in \cite[Thm. 8.18]{folland1999real}. Therefore, we can replace function $f$ in \eqref{ocp:dyn} with $\omega f$, which is compactly supported in $\bRnx$, ensuring Assumption \ref{asm:x-bnd}. Note that the inclusion of $\omega$ is only for the purpose of proofs of results discussed in subsequent sections; it is not necessary in a practical implementation.

The optimal control problem \eqref{ocp} is posed in the Mayer form \cite[Sec. 3.3.2]{liberzon2011calculus}, without loss of generality. The original problem could be provided in the Bolza or Lagrange forms where the objective function contains the integral of a continuously differentiable running cost function $L_{\mathrm{r}}:\bRp\times\bR^{n_x}\times\bR^{n_u}\to\bR$, given by
$$
    \int_{\tinit}^{\tfinal}L_{\text{r}}(t,x(t),u(t))\mathrm{d}t
$$
which is well-defined when $x$ is a state trajectory and $u$ is the corresponding piecewise continuous control input. The Mayer form is obtained from the Lagrange and Bolza forms by equivalently reformulating the integral of $L_{\mathrm{r}}$ as a new dynamical system
\begin{equation}
    \dot{l}(t) = L_{\mathrm{r}}(t,x(t),u(t))\label{run-cost-dyn}
\end{equation}
and adding $l(\tfinal)$ as a terminal state cost to the objective function, along with specifying the initial condition $l(\tinit) = 0$.

We propose an SCP and multiple-shooting based solution method \cite{bock1984multiple} (which only requires first-order information). The gradient of the integral of $L_{\mathrm{r}}$ can be computed through the first-order sensitivities \cite[Sec. 4.2]{quirynen2015autogenerating} of dynamical system \eqref{run-cost-dyn}---a qualitatively similar effort to linearizing \eqref{ocp:dyn}. Therefore, converting to the Mayer form is beneficial for consolidating majority of the effort of gradient computation in one place, which is efficient for numerical implementations. Furthermore, the integral of running cost function can be either convex or nonconvex. In the convex case, however, the integral often cannot be reformulated as a canonical conic constraint with the help of additional slack variables (see Section \ref{subsec:cvx-exact-penalty} for further details). Moreover, a convex running cost function in a free-final-time problem will become nonconvex after conversion to a fixed-final-time problem via the time-dilation technique described in Section \ref{subsec:dilation}. In the subsequent development, unless stated otherwise, we assume that the dynamical system \eqref{ocp:dyn}, the terminal state cost function $L$, and the boundary condition \eqref{ocp:bc-eq} already encode any desired running cost terms. 
%
\subsection{Constraint Reformulation}\label{subsec:cnstr-reform}
To reformulate path constraints \eqref{ocp:ineq-cnstr} and \eqref{ocp:eq-cnstr}, we consider continuous \textit{exterior penalty functions}
\begin{subequations}
    \begin{align}
   &  q_i : \bRp\times\bRnx\times\bRnu \to \bRp,~~i=1,\ldots,\Nineq \\ 
   & q_i(z)  \begin{cases}
                    {}= 0 &\text{if }z \le 0 \\
                    {}> 0 & \text{otherwise}
                \end{cases} \\
    & p_j : \bRp\times\bRnx\times\bRnu \to \bRp,~~j=1,\ldots,\Neq \\ 
    & p_j(z)  \begin{cases}
                    {}= 0 &\text{if }z = 0 \\
                    {}> 0 & \text{otherwise}
                \end{cases}
    \end{align}\label{atom-ext-penal}%
\end{subequations}
and the penalty function $\penaltyfun : \bRp\times\bRnx\times\bRnu\to\bRp$ 
\begin{align}
  & \penaltyfun(t,x,u) = \sum_{i=1}^{\Nineq} q_i(g_i(t,x,u)) + \sum_{j=1}^{\Neq} p_j(h_j(t,x,u))\label{P-def}
\end{align}
Note that $g_i$, $i=1,\ldots,\Nineq$, and $h_j$, $j=1,\ldots,\Neq$, denote the scalar-valued elements of functions $g$ and $h$, respectively.
It is straightforward to show that $t\mapsto \penaltyfun(t,x(t),u(t))$ is piecewise continuous over $\tspan$ when $x$ is a state trajectory and $u$ is the corresponding piecewise continuous control input.
\begin{lem}\label{lem:cnstr-intgl} Path constraints \eqref{ocp:ineq-cnstr} and \eqref{ocp:eq-cnstr} are satisfied by $x(t)$ and $u(t)~\ae~t\in\tspan$ iff
\begin{align}
 & \int_{\tinit}^{\tfinal}\penaltyfun(t,x(t),u(t))\mathrm{d}t = 0 \label{cnstr-intgl}
\end{align}
\end{lem}\vspace{-0.5cm}
\begin{pf}
For each $t\in\tspan$, $\penaltyfun(t,x(t),u(t))$ is nonnegative since $q_i(g_i(t,x(t),u(t)))$ and $p_j(h_j(t,x(t),u(t)))$ are nonnegative. Since $t\mapsto \penaltyfun(t,x(t),u(t))$ is piecewise continuous, \eqref{cnstr-intgl} holds iff $\penaltyfun(t,x(t),u(t)) = 0$ $\ae$ $t\in\tspan$.\\
Then, $q_i(g_i(t,x(t),u(t))) = 0$ and $p_j(h_j(t,x(t),u(t))) = 0$ $\ae$ $t\in\tspan$, which is equivalent to $g_i(t,x(t),u(t)) \le 0$ and $h_j(t,x(t),u(t)) = 0$ $\ae$ $t\in\tspan$.
\end{pf}
Next, we augment the system dynamics with an additional state to accumulate the constraint violations, and impose periodic boundary conditions on it to ensure continuous-time constraint satisfaction.
\begin{cor}
The differential-algebraic system
\begin{subequations}
\begin{align}
 & \dot{x}(t) = f(t,x(t),u(t))\\
 & g(t,x(t),u(t)) \le 0\label{cnstr-reform:ineq}\\
 & h(t,x(t),u(t)) = 0\label{cnstr-reform:eq}
\end{align}\label{cnstr-reform}%
\end{subequations}
is satisfied by $x$ and $u$ a.e. in $\tspan$ iff $x$ and $u$ solve the boundary value problem on $\tspan$
\begin{subequations}
\begin{align}
\dot{x}(t) ={} & f(t,x(t),u(t))\label{tpbvp:original-dyn}\\
\dot{y}(t) ={} & \penaltyfun(t,x(t),u(t))\label{tpbvp:aug-dyn}\\
y(0) ={} & y(\tfinal)\label{tpbvp:bc}
\end{align}\label{tpbvp}%
\end{subequations}
\end{cor}
\begin{pf} Suppose that $x$ and $u$ satisfy \eqref{ocp:dyn}, \eqref{ocp:ineq-cnstr}, and \eqref{ocp:eq-cnstr} almost everywhere in  $\tspan$. Since $t\mapsto \penaltyfun(t,x(t),u(t))$ is piecewise continuous over $\tspan$, there exists a unique state trajectory $y:\tspan\to\bR$ which satisfies \eqref{tpbvp:aug-dyn} almost everywhere in $\tspan$ with initial condition $y(\tinit) = 0$. Then, $y(\tfinal) = 0$ iff \eqref{cnstr-intgl} holds. Therefore, from Lemma \eqref{lem:cnstr-intgl}, the path constraints \eqref{ocp:ineq-cnstr} and \eqref{ocp:eq-cnstr} are satisfied almost everywhere in $\tspan$.\\
For the other direction, we let $x$, $y$, and $u$ be the solution to \eqref{tpbvp}. Then, \eqref{tpbvp:bc} implies that \eqref{cnstr-intgl} holds. We obtain the desired result from Lemma \ref{lem:cnstr-intgl}.
\end{pf}
\begin{rem}\label{rem:dyn-jac-cont} Since  a gradient-based solution method will be adopted,  the exterior penalty functions $q_i$ and $p_j$, for $i=1,\ldots,\Nineq$ and $j=1,\dots,\Neq$, must be  continuously differentiable in order to compute first-order sensitivities \cite[Sec. 4.2]{quirynen2015autogenerating} of \eqref{tpbvp:aug-dyn}. For e.g., 
%
\begin{align}
    q_i(z) ={} & \normplus{z}^2,~p_j(z) = z^2\label{choose-ext-penal}
\end{align}   
for any $z\in\bR$, are valid choices. Then, $\penaltyfun$ can be compactly written as
\begin{align}
    \penaltyfun(t,x,u) = \ones{n_g}^\top\normplus{g(t,x,u)}^2 + \ones{n_h}^\top h(t,x,u)^2\label{P-choose-ext-penal}
\end{align}
for any $t\in\bR$, $x\in\bRnx$, and $u\in\bRnu$.
\end{rem}
\begin{rem}\label{rem:vanish-grad-P}
The partial derivatives of $\penaltyfun$ are given by
\begin{align}
    \inlpbyp{\penaltyfun}{\square} = \sum_{i=1}^{\Nineq}\nabla q_i\inlpbyp{g_i}{\square} + \sum_{i=1}^{\Nineq}\nabla p_j\inlpbyp{h_j}{\square}\label{P-def-partial}
\end{align}
where $\square = t$, $x$, and $u$. A consequence of the continuous differentiability of $q_i$ and $p_j$, and their construction in \eqref{atom-ext-penal}, is that their derivatives in \eqref{P-def-partial} are zero wherever $g_i \le 0$ and $h_j = 0$.     
\end{rem}
\begin{rem}
The penalty function $\penaltyfun$ can be defined more generally as a vector-valued function of the form
\begin{align}
    \penaltyfun(t,x,u) = M\left[\begin{array}{c}q(g(t,x,u))\\p(h(t,x,u))\end{array}\right]\label{P-def-gen}
\end{align}
for any $t\in\bRp$, $x\in\bRnx$, and $u\in\bRnu$, where matrix $M \in \bR^{n_y \times (\Nineq+\Neq)}$, which we refer to as the mixing matrix, has nonnegative entries, and functions $q : \bR^{\Nineq} \to \bR^{\Nineq}$ and $p:\bR^{\Neq}\to \bR^{\Neq}$ are defined as
\begin{subequations}
\begin{align}
    q((z_g^1,\ldots,z_g^{\Nineq})) &{}= (q_1(z_g^1),\ldots,q_{\Nineq}(z_g^{\Nineq}))\\
    p((z_h^1,\ldots,z_h^{\Neq})) &{}= (p_1(z_h^1),\ldots,p_{\Neq}(z_h^{\Neq}))
\end{align}%
\end{subequations}
for any $z^i_g \in \bR$ and $z_h^j \in \bR$, with $i=1,\ldots,\Nineq$ and $j=1,\ldots,\Neq$. The mixing matrix has at most $\Nineq+\Neq$ rows, at most one positive entry per column, and exactly $\Nineq+\Neq$ positive entries in total. We obtain \eqref{P-def} and a scalar-valued state $y$ in \eqref{tpbvp:aug-dyn} by choosing $M = \ones{\Nineq+\Neq}^\top$. In general, choosing \eqref{P-def-gen} will lead to a vector-valued state $y\in\bR^{n_y}$ for quantifying the total constraint violation. While \eqref{P-def} is suitable for most cases, choosing \eqref{P-def-gen} for certain applications can be numerically more reliable. In the remainder of the discussion, we adopt the definition in \eqref{P-def}. However, the results also apply to the general representation in \eqref{P-def-gen}.
\end{rem}
%
\subsection{Generalized Time-Dilation}\label{subsec:dilation}
%
\textit{Time-dilation} is a transformation technique for posing a free-final-time optimal control problem as an equivalent fixed-final-time problem \cite[Sec. III.A.1]{szmuk2020successive}. Let $x$ be a state trajectory for \eqref{ocp:dyn} over $\tspan$ obtained with control input $u$ and some initial condition. We define a strictly increasing, continuously differentiable mapping $t : [0,1]\to\bRp$, with boundary conditions $t(0)=\tinit,\,t(1)=\tfinal$, and treat the derivative of the map 
$$
    s(\tau) = \dbyd{t(\tau)}{\tau} = \derv{t}(\tau) 
$$
for $\tau\in[0,1]$, as an additional control input, which we refer to as the \textit{dilation factor}. Note that the derivative with respective $\tau\in[0,1]$ is denoted by $\derv{\square} = \dbyd{\square}{\tau}$. Previously, dilation factors were either treated as a single parameter decision variable (as in \cite[Sec. III.A.1]{szmuk2020successive}) or multiple discrete parameter decision variables (as in \cite[Sec. 2.1]{kamath2023seco}). The approach proposed herein generalizes it by treating dilation factor as a continuous-time control input. We note that this approach bears resemblance to the time-scaling transformation proposed in \cite{lee1997control,lin2014control}. Further, we treat $t$ as an additional state, and define functions $\tilde{u}$ and $\tilde{x}$ as 
\begin{subequations}
\begin{align}
& \tilde{u}(\tau) = \big(u(t(\tau)),\,s(\tau)\big)\in\bR^{n_u+1}\\
& \tilde{x}(\tau) = \big(x(t(\tau)),\,y(t(\tau)),\,t(\tau)\big)\in\bR^{n_x+2}
\end{align}
\end{subequations}
for $\tau\in[0,1]$. Then, using chain-rule, we obtain 
\begin{align}
\derv{\tilde{x}}(\tau) = \begin{bmatrix} \dot{x}(t) \\ \dot{y}(t) \\ 1 \end{bmatrix} \dbyd{t(\tau)}{\tau} = F(\tilde{x}(\tau),\tilde{u}(\tau))\label{aug-dyn}
\end{align}
for $\tau\in[0,1]$, where 
\begin{align}
& F(\tilde{x}(\tau),\tilde{u}(\tau)) = 
\begin{bmatrix}f(t(\tau),x(t(\tau)),u(t(\tau)))\\\penaltyfun(t(\tau),x(t(\tau)),u(t(\tau)))\\1\end{bmatrix}\!s(\tau)\label{aug-dyn-RHS}
\end{align}
We refer to \eqref{aug-dyn} as the \textit{augmented system}, with \textit{augmented state} $\tilde{x}$, and \textit{augmented control input} $\tilde{u}$.

In the subsequent development, we define $\ntildx=n_x+2$ and use matrices $\selector{x}$, $\selector{y}$, and $\selector{t}$ to select the first $n_x$ elements ($x$), the penultimate element ($y$), and the last element ($t$), respectively, of $\tilde{x}$. Similarly, matrices $\selector{u}$ and $\selector{s}$ select the first $n_u$ elements ($u$) and the last element ($s$), respectively, of $\tilde{u}$. Note that generalized time-dilation transforms non-autonomous dynamical systems into equivalently autonomous ones (without explicit dependence on time) by increasing the dimensions of the state and control input. If the final time is fixed in \eqref{ocp}, then the generalized time-dilation step can be skipped. Henceforth, we drop the qualifier ``generalized'' for brevity.
\begin{rem}
When the system dynamics and constraint functions in \eqref{ocp} are time-varying, we use time-dilation to augment the system with an extra state and control input. However, when the system and constraint descriptions are time-invariant, it suffices to augment only the control input.     
\end{rem} 
%
\subsection{Reformulated Optimal Control Problem}
%
%
Subjecting \eqref{ocp} to constraint reformulation and time-dilation results in the following optimal control problem
\begin{eqbox}
\begin{subequations}
\begin{align}
\underset{\tilde{x},\,\tilde{u}}{\operatorname{minimize}}&~\tilde{L}(\tilde{x}(1)) \\
\operatorname{subject~to}&~\derv{\tilde{x}}(\tau) = F(\tilde{x}(\tau),\tilde{u}(\tau)) & & \!\!\!\!\ae~\tau\in[0,1]\label{ocp-reform:dyn}\\
 &~\selector{s}\tilde{u}(\tau)>0 & & \!\!\!\!\phantom{\ae~}\tau\in[0,1] \label{ocp-reform:bound-dilation}\\
 &~\selector{y}(\tilde{x}(1)-\tilde{x}(0)) = 0\label{ocp-reform:cnstr-bc}\\
 &~\tilde{P}(\tilde{x}(0),\tilde{x}(1)) \le 0 \label{ocp-reform:bc-ineq}\\
 &~\tilde{Q}(\tinit,\tilde{x}(0),\tilde{x}(1)) = 0 \label{ocp-reform:bc-eq}
\end{align}\label{ocp-reform}%
\end{subequations}
\end{eqbox}
where, 
\begin{subequations}
\begin{align}
\tilde{L}(\tilde{x}) ={} & L(\selector{t}\tilde{x},\selector{x}\tilde{x})\\
\tilde{P}(\tilde{x},\tilde{x}^\prime) ={} & P(\selector{t}\tilde{x},\selector{x}\tilde{x},\selector{t}\tilde{x}^\prime,\selector{x}\tilde{x}^\prime)\\
\tilde{Q}(\tinit,\tilde{x},\tilde{x}^\prime) ={} & \begin{bmatrix}Q(\tinit,\selector{x}\tilde{x},\selector{t}\tilde{x}^\prime,\selector{x}\tilde{x}^\prime)\\\selector{t}\tilde{x} - \tinit\end{bmatrix}
\end{align}
\end{subequations}
for any $\tilde{x},\tilde{x}^\prime \in \mathbb{R}^{n_{\tilde{x}}}$. Since time is treated as a state variable, its initial condition needs to be specified through $\tilde{Q}$. This is not needed in \eqref{ocp} because (constant) $\tinit$ is explicitly passed to functions $P$ and $Q$. Furthermore, positivity of the dilation factor is ensured using \eqref{ocp-reform:bound-dilation}.
Observe that, as a consequence of the constraint reformulation, the continuous-time constraints on the state and control input are not explicitly imposed in \eqref{ocp-reform}; they are instead embedded within the augmented system in \eqref{ocp-reform:dyn} and boundary condition \eqref{ocp-reform:cnstr-bc}.
%
\section{Parameterization and Discretization}
This section describes the transformation of the (infinite-dimensional) reformulated optimal control problem \eqref{ocp-reform} into a (finite-dimensional) nonconvex optimization problem through \textit{parameterization} of the augmented control input, and \textit{discretization} of \eqref{ocp-reform:dyn} over the interval $[0,1]$. In addition, we introduce a relaxation to the constraint reformulation step after discretization to ensure favorable properties from an optimization viewpoint, and analyze the consequences of the relaxation.

We parameterize the augmented control input with finite-dimensional vectors via $\tilde{\nu}:[0,1]\to\bR^{n_{\tilde{u}}}$, defined as
\begin{equation}
\tilde{\nu}(\tau) = \begin{pmatrix}\Gamma_u(\tau)\otimes\begin{bmatrix}
     \eye{n_u}  \\
     \zeros{1\times n_u}
\end{bmatrix}\end{pmatrix} U + \begin{bmatrix}
     \zeros{n_u\times N_s}  \\
     \Gamma_s(\tau)  
\end{bmatrix}S\label{aug-ctrl-param}
\end{equation}
for $\tau\in[0,1]$, where 
\begin{equation}
    \Gamma_{\scriptscriptstyle\square} = \big[\Gamma_{\scriptscriptstyle\square}^1\,\ldots\,\Gamma_{\scriptscriptstyle\square}^{N_{\scriptscriptstyle\square}}\big]\label{basis-func}
\end{equation}
Each element in \eqref{basis-func}, $\Gamma_{\scriptscriptstyle\square}^k : [0,1] \to \bR$, is a polynomial basis function, for $k=1,\ldots,N_{\scriptscriptstyle\square}$ and $\square = u,s$. Vectors $U\in\bR^{n_uN_u}$ and $S\in\bR^{N_s}$, which are the decision variables, are coefficients parameterizing the control input and dilation factor. The support of each of the basis functions in $\Gamma_u$ and $\Gamma_s$ can be pairwise disjoint and strictly contained in $[0,1]$. 
\begin{equation}
    \tilde{\nu}(\tau) = \big( (\Gamma_u(\tau) \otimes \eye{n_u})U,\,\Gamma_s(\tau)S\big)\label{aug-ctrl-param-compact}
\end{equation}
\begin{rem}
The parameterization \eqref{aug-ctrl-param} for the augmented control input is general enough to subsume several control input parameterizations used in practice, such as zero- and first-order-hold, and pseudospectral polynomials in orthogonal collocation (see Appendix \ref{ctrl-param-eg}). We distinguish between the parameterizations for the control input and dilation factor because the former is allowed to have a general polynomial parameterization as long as $U$ is elementwise bounded, whereas, the latter must be constrained to be bounded and positive due to \eqref{ocp-reform:bound-dilation}. Such a requirement should be straightforward to handle with convex constraints on $S$.
\end{rem}
After parameterizing the augmented control input, we discretize \eqref{ocp-reform:dyn} via its integral form over a finite grid of size $N$ in $[0,1]$, denoted by $0=\tau_1<\ldots<\tau_N=1$. We treat the augmented states $\tilde{x}_k$ at the node points $\tau_k$ also as decision variables, and denote them compactly as 
$$
    \tilde{X} = (\tilde{x}_1,\ldots,\tilde{x}_N)\in\bR^{\ntildx N}
$$
Next, we define $F_k:\bR^{\ntildx} \times \bR^{\ntildx} \times \bR^{n_uN_u} \times \bR^{N_s} \to \bR^{\ntildx}$ 
\begin{align}
 & (\tilde{x}_{k+1},\tilde{x}_{k},U,S) \mapsto \label{Fk}\\
 & ~~\qquad\qquad\tilde{x}_{k+1}-\tilde{x}_{k}-\int_{\tau_k}^{\tau_{k+1}}F\big(\tilde{x}^k(\tau),\tilde{\nu}(\tau)\big)\mathrm{d}\tau \nonumber 
\end{align}
for $k=1,\ldots,N-1$, where the augmented state trajectory $\tilde{x}^{k}$ satisfies \eqref{ocp-reform:dyn} almost everywhere on $[\tau_k,\tau_{k+1}]$ with augmented control input $\tilde{\nu}$, and initial condition $\tilde{x}_k$. Then, the discretization of \eqref{ocp-reform:dyn}, together with \eqref{ocp-reform:cnstr-bc}, yields the constraints
\begin{subequations}
\begin{align}
& F_k(\tilde{x}_{k+1},\tilde{x}_k,U,S) = 0\label{disc-dyn}\\
& \selector{y}(\tilde{x}_{k+1}-\tilde{x}_{k})=0\label{disc-cnstr-bc}
\end{align}\label{disc}%
\end{subequations}
for $k=1,\ldots,N-1$. The discretization in \eqref{disc-dyn}, commonly referred to as \textit{multiple-shooting}, is exact because there is no approximation involved in switching to the integral representation of the differential equation \eqref{ocp-reform:dyn}, and it is numerically more stable compared to its precursor, single-shooting \cite[Sec. VI.A.3]{betts1998survey}.  
%
\begin{defn}[Continuous-Time Feasibility]\label{def:ct-feas} The triplet $\tilde{X}$, $U$, and $S$ is said to be continuous-time feasible if 
it satisfies \eqref{disc}.
\end{defn}
A continuous-time feasible triplet: $\tilde{X}$, $U$, and $S$ corresponds to an augmented state trajectory generated by  integration of $\eqref{ocp-reform:dyn}$ over $[0,1]$ with augmented control input $\tilde{\nu}$ and initial condition $\tilde{x}_1$, such that it passes through the node values $\tilde{x}_k$ at $\tau_k$, for $k=2,\ldots,N$. Furthermore, due to \eqref{disc-cnstr-bc}, the resulting state trajectory and control input satisfy the path constraints almost everywhere in $\tspan$.

Continuous-time feasibility is a crucial metric for assessing the quality of solutions computed by direct methods. They solve for finite-dimensional quantities (such as $\tilde{X}$, $\tilde{U}$, and $S$) which parameterize the control input, and at times, even the state trajectory. In particular, direct collocation methods \cite{rao2009survey} approximate the state trajectory with a basis of polynomials (for e.g. orthogonal collocation methods \cite{garg2010unified,ross2012review}). They match the derivative of the polynomial approximant to the right-hand-side of \eqref{ocp:dyn}, and impose the path constraints at finitely-many nodes. As a result, such methods cannot guarantee that the computed solution is continuous-time feasible, necessitating computationally intensive heuristics such as mesh refinement \cite[Sec. 4.7]{betts2010practical}. On the other hand, the proposed framework leverages multiple-shooting and constraint reformulation to provide continuous-time feasible solutions, which is a key distinction  compared to existing direct methods \cite[Sec. 2.4.1]{malyuta2021advances}, \cite[Table 1.1]{chai2023advanced}.
%
\subsection{Constraint Qualification and Relaxation}\label{subsec:cnstr-qual-relax}
An issue with directly imposing \eqref{disc} is that it violates \textit{linear independence constraint qualification (LICQ)} \cite[Def. 12.4]{nocedal2006numerical}. For each point in $\bR^{n_{\tilde{x}}N}\times\bR^{n_uN_u}\times\bR^{N_s}$, we can determine the active set \cite[Def. 12.1]{nocedal2006numerical} corresponding to the constraints in \eqref{disc}. LICQ is said to hold at a point if the gradients of the constraints in the corresponding active set are linearly independent. The following result shows that LICQ is violated at all points feasible with respect to \eqref{disc}.
\begin{lem} If $\tilde{X}$, $U$, and $S$ are feasible with respect to \eqref{disc-dyn} and \eqref{disc-cnstr-bc}, then for $k=1,\ldots,N-1$
\begin{subequations}
\begin{align}
\selector{y}\inlpbyp{F_k}{\tilde{x}_{k+1}}(\tilde{x}_{k+1},\tilde{x}_{k},U,S) ={} & [\zeros{1\times n_x}~1~0]\label{jac-Fk-nonzero-xk}\\
\selector{y}\inlpbyp{F_k}{\tilde{x}_k}(\tilde{x}_{k+1},\tilde{x}_{k},U,S) ={} & [\zeros{1\times n_x}-1~0]\label{jac-Fk-nonzero-xkp1}\\
 \selector{y}\inlpbyp{F_k}{\scriptscriptstyle U}(\tilde{x}_{k+1},\tilde{x}_{k},U,S) ={} & \zeros{1\times n_uN_u}\\
 \selector{y}\inlpbyp{F_k}{\scriptscriptstyle S}(\tilde{x}_{k+1},\tilde{x}_{k},U,S) ={} & \zeros{1\times N_s}
\end{align}\label{jac-Fk-nonzero}%
\end{subequations}\label{lem:cnstr-dyn-grad}%
\end{lem}
\begin{pf}
Suppose $\tilde{X}$, $U$, and $S$ satisfy \eqref{disc-dyn} and \eqref{disc-cnstr-bc}. For each $k=1,\ldots,N-1$, let $\tilde{x}^k$ be the augmented state trajectory on $[\tau_k,\tau_{k+1}]$ with augmented control input $\tilde{\nu}$ and initial condition $\tilde{x}_k$. Then, due to Lemma \ref{lem:cnstr-intgl}, path constraints \eqref{ocp:ineq-cnstr} and \eqref{ocp:eq-cnstr} are satisfied almost everywhere on the time interval $[\selector{t}\tilde{x}_k,\selector{t}\tilde{x}_{k+1}]$. Recall that $F$ defined in \eqref{aug-dyn-RHS} is continuously differentiable and that $\tilde{\nu}$ is piecewise continuous. This  implies that 
\begin{align*}
 & \tau\mapsto\inlpbyp{F(\tilde{x}^k(\tau),\tilde{\nu}(\tau))}{\tilde{x}},~\tau\mapsto\inlpbyp{F(\tilde{x}^k(\tau),\tilde{\nu}(\tau))}{\tilde{u}}
\end{align*}
are piecewise continuous functions. Therefore, they are integrable and the partial derivatives of $F_k$ are well-defined. Elements of the penultimate row of $\inlpbyp{F}{\tilde{x}}$ and $\inlpbyp{F}{\tilde{u}}$ 
(i.e., the partial derivatives of $\penaltyfun$), are zero $\ae$ on $[\tau_k,\tau_{k+1}]$ when evaluated on  $\tilde{x}^k$ and $\tilde{\nu}$. This is because all elements of this row contain factors of the derivatives of $q_i$ or $p_j$, which are zeros due to Remark \ref{rem:vanish-grad-P}. Hence, only the partial derivatives of $F_k$ with respect to $\tilde{x}_{k}$ and $\tilde{x}_{k+1}$ contain nonzero elements in that row, as shown in \eqref{jac-Fk-nonzero-xk} and \eqref{jac-Fk-nonzero-xkp1}.
\end{pf}
Using Lemma \ref{lem:cnstr-dyn-grad}, for  any $\tilde{X}$, $U$, and $S$ satisfying \eqref{disc-dyn} and \eqref{disc-cnstr-bc}, the gradient of the left-hand-side of \eqref{disc-cnstr-bc}, and the penultimate row of the gradient of $F_k$ (shown in \eqref{jac-Fk-nonzero}) are identical. Hence, all feasible solutions of an optimization problem having both \eqref{disc-dyn} and \eqref{disc-cnstr-bc} as constraints will not satisfy LICQ. 
The proposed framework relies on the classical exact penalization result to determine KKT points of \eqref{ocp-reform-ctrlprm} \cite[Thm. 17.4]{nocedal2006numerical}, which are related to local minimizers when LICQ holds \cite[Thm. 12.1]{nocedal2006numerical}. So, we remedy the pathological scenario due to \eqref{disc} by relaxing \eqref{disc-cnstr-bc} to an inequality with a positive constant $\epsilon$. 

Then, \eqref{ocp-reform} transforms under augmented control input parameterization \eqref{aug-ctrl-param} and discretization \eqref{disc} as follows 
\begin{eqbox}
\begin{subequations}
\begin{align}
\underset{\tilde{X},\,U,\,S}{\operatorname{minimize}}~~&~\tilde{L}(\tilde{x}_N)\\
\operatorname{subject~to}~&~F_k(\tilde{x}_{k+1},\tilde{x}_k,U,S) = 0\label{ocp-reform-ctrlprm:dyn} \\
 &~\selector{y}(\tilde{x}_{k+1}-\tilde{x}_k) \le \epsilon \label{ocp-reform-ctrlprm:relax}\\
 &~k=1,\ldots,N-1 \nonumber\\
 &~U\in\mc{U},~S\in\mc{S}\label{ocp-reform-ctrlprm:ctrlcvx}\\
 &~\tilde{P}(\tilde{x}_1,\tilde{x}_N) \le 0,~\tilde{Q}(\tinit,\tilde{x}_1,\tilde{x}_N) = 0\label{ocp-reform-ctrlprm:bc}
\end{align}\label{ocp-reform-ctrlprm}%
\end{subequations}
\end{eqbox}
where $\mc{U}\subset\bR^{n_uN_u}$ and $\mc{S}\subset\bR^{N_s}$ are compact convex sets. 
\begin{rem}
If $\tilde{X}$, $U$, and $S$ satisfy \eqref{ocp-reform-ctrlprm:dyn} and \eqref{ocp-reform-ctrlprm:relax}, then LICQ will not be trivially violated by the active set associated with \eqref{ocp-reform-ctrlprm:relax} and the penultimate row of \eqref{ocp-reform-ctrlprm:dyn}. If \eqref{ocp-reform-ctrlprm:relax} holds with strict inequality, for all $k=1,\ldots,N-1$, then the result follows immediately. However, if for some $k=1,\ldots,N-1$, \eqref{ocp-reform-ctrlprm:relax} holds with equality, then
\begin{equation}
    \int_{\tau_k}^{\tau_{k+1}}\!\selector{s}\tilde{\nu}(\tau)\penaltyfun\big(\selector{t}\tilde{x}(\tau),\selector{x}\tilde{x}^k(\tau),\selector{u}\tilde{\nu}(\tau)\big)\mathrm{d}\tau = \epsilon\label{ocp-reform-ctrlprm:relax-active}
\end{equation}
where $\tilde{x}^k$ is the augmented state trajectory with augmented control input $\tilde{\nu}$ and initial condition $\tilde{x}_k$. Due to piecewise continuity of the integrand in \eqref{ocp-reform-ctrlprm:relax-active} with respect to $\tau\in[\tau_k,\tau_{k+1}]$, it implies that there exists an interval $\mc{I}_k \subset [\tau_k,\tau_{k+1}]$ where the integrand of \eqref{ocp-reform-ctrlprm:relax-active} is positive, i.e., some of the path constraints are violated on the segments of $\tilde{x}^k$ and $\tilde{\nu}$ corresponding to interval $\mc{I}_k$. Therefore, unlike the case in Lemma \ref{lem:cnstr-dyn-grad}, the partial derivatives of $\penaltyfun$ evaluated on those segments of $\tilde{x}^k$ and $\tilde{\nu}$ will not be trivially zero. 
\end{rem}
\begin{rem}
Besides addressing the LICQ issue, relaxation \eqref{ocp-reform-ctrlprm:relax} is essential for enabling exact penalization---a key step in the proposed framework. The combination of \eqref{disc-dyn} and \eqref{disc-cnstr-bc} is equivalent to the following constraint for each $k=1,\ldots,N-1$
\begin{align}
    \selector{y}(\tilde{x}_{k+1}-\tilde{x}_k-F_k(\tilde{x}_{k+1},\tilde{x}_k,U,S)) = 0\label{ocp-reform-ctrlprm:unrelax}
\end{align}
Exact penalization of \eqref{ocp-reform-ctrlprm:unrelax} involves penalizing the left-hand-side of \eqref{ocp-reform-ctrlprm:unrelax} using a nonsmooth exterior penalty function (such as $p_j$). However, the effect of exterior penalty is suppressed due to the positivity and continuous differentiability of $\penaltyfun$, which destroys the exactness of the penalty term (see discussion in \cite[p. 513]{nocedal2006numerical}). In other words, it may not be possible to recover a solution that is feasible with respect to \eqref{ocp-reform-ctrlprm:unrelax} using a finite weight for the penalty term.
\end{rem}
\begin{rem}\label{rem:UScnstr}
Set $\mc{U}$ ensures that the control input lies in a compact convex set. Whenever the choice of parameterization permits (e.g., with zero- or first-order-hold), the convex control constraints in \eqref{ocp:ineq-cnstr} and \eqref{ocp:eq-cnstr} could be captured with $\mc{U}$, i.e., via imposing convex control constraints only at the node points (inter-sample constraint satisfaction is guaranteed); otherwise, $\mc{U}$ will encode elementwise bounds on $U$. Set $\mc{S}$ ensures that the dilation factor is positive and bounded. Besides handling convex control constraints directly, requiring the augmented control input to lie in a compact set via \eqref{ocp-reform-ctrlprm:ctrlcvx} is necessary for deriving a pointwise constraint violation bound, given the $\epsilon$-relaxation in \eqref{ocp-reform-ctrlprm:relax}. 
\end{rem}
%
\subsection{Gradient of Discretized Dynamics}\label{subsec:jacob-dyn}
A key step in any gradient-based solution method for \eqref{ocp-reform-ctrlprm} is to compute the partial derivatives of $F_k$. We adopt the so-called variational method \cite[Sec. 3.2]{lin2014control}, \cite[Sec. 4.2]{quirynen2015autogenerating}, also referred to as inverse-free exact discretization \cite[Sec. 2.3]{kamath2023seco}. 

Let $\bar{\tilde{X}} = (\bar{\tilde{x}}_1,\ldots,\bar{\tilde{x}}_N)$, $\bar{U}$, and $\bar{S}$ denote a reference solution, and let $\bar{\tilde{\nu}}$ denote the corresponding parameterization using \eqref{aug-ctrl-param}. For each $k=1,\ldots,N-1$, consider the following initial value problem over $[\tau_k,\tau_{k+1}]$
\begin{subequations}
\begin{align}
 & \derv{\Phi}{}^k_{\tilde{x}}(\tau) = A(\tau)\Phi^k_{\tilde{x}}(\tau)\\[-0.1cm]
 & \derv{\Phi}{}^k_{u}(\tau) = A(\tau)\Phi^k_{u}(\tau) + \Gamma_u(\tau)\otimes B(\tau)\\[-0.1cm]
 & \derv{\Phi}{}^k_{s}(\tau) = A(\tau)\Phi^k_{s}(\tau) + \Gamma_s(\tau)\otimes C(\tau)\\[-0.1cm]
 & \Phi_{\tilde{x}}^k(\tau_k) = \eye{\ntildx}\phantom{\derv{\Phi}}\\[-0.1cm]
 & \Phi_u^k(\tau_k) = \zeros{\ntildx\times n_uN_u}\phantom{\derv{\Phi}}\\[-0.1cm]
 & \Phi_s^k(\tau_k) = \zeros{\ntildx\times N_s}\phantom{\derv{\Phi}}
\end{align}\label{sensitivity-ivp}%
\end{subequations}
where $\Phi_{\tilde{x}}^k(\tau)\in\bR^{\ntildx\times\ntildx}$, $\Phi_u^k(\tau)\in\bR^{\ntildx\times n_uN_u}$, and $\Phi_s^k(\tau)\in\bR^{\ntildx\times N_s}$, for $\tau\in[\tau_k,\tau_{k+1}]$. The partial derivatives of $F$ evaluated on the augmented state trajectory $\bar{\tilde{x}}^k$ for \eqref{ocp-reform:dyn} over $[\tau_k,\tau_{k+1}]$, with augmented control input $\bar{\tilde{\nu}}$ and initial condition $\bar{\tilde{x}}_k$, are compactly denoted with $A(\tau)\in\bR^{\ntildx\times\ntildx}$, $B(\tau)\in\bR^{\ntildx\times n_u}$, and $C(\tau)\in\bR^{\ntildx}$, i.e.,
\begin{subequations}
\begin{align}
 A(\tau) ={} & \inlpbyp{F}{\tilde{x}}\big(\bar{\tilde{x}}^k(\tau),\bar{\tilde{\nu}}(\tau)\big)\\
 B(\tau) ={} & \inlpbyp{F}{\tilde{u}}\big(\bar{\tilde{x}}^k(\tau),\bar{\tilde{\nu}}(\tau)\big)\selector{u}^\top\\
 C(\tau) ={} & \inlpbyp{F}{\tilde{u}}\big(\bar{\tilde{x}}^k(\tau),\bar{\tilde{\nu}}(\tau)\big)\selector{s}^\top
\end{align}%
\end{subequations}
for $\tau\in[\tau_k,\tau_{k+1}]$. Then, the partial derivatives of $F_k$ with respect to $\tilde{x}_k$, $U$, and $S$ (denoted by $A_k$, $B_k$, and $C_k$) evaluated at the reference solution are related to the terminal value of the solution to \eqref{sensitivity-ivp} as follows
\begin{subequations}
\begin{align}
A_k ={} & \Phi^k_{\tilde{x}}(\tau_{k+1}) = - \inlpbyp{F_k}{\tilde{x}_{k}}(\bar{\tilde{x}}_{k+1},\bar{\tilde{x}}_k,\bar{U},\bar{S}) \\
B_k ={} & \Phi_u^k(\tau_{k+1}) = -\inlpbyp{F_k}{\scriptscriptstyle U}(\bar{\tilde{x}}_{k+1},\bar{\tilde{x}}_k,\bar{U},\bar{S})\\
C_k ={} & \Phi_{s}^k(\tau_{k+1}) = -\inlpbyp{F_k}{\scriptscriptstyle S}(\bar{\tilde{x}}_{k+1},\bar{\tilde{x}}_k,\bar{U},\bar{S})
\end{align}
\end{subequations}
Finally, the linearization of $F_k$ is given by the mapping
\begin{align}
& \hspace{-0.5cm}(\tilde{x}_{k+1},\tilde{x}_{k},U,S) \mapsto \tilde{x}_{k+1} - \bar{\tilde{x}}^k(\tau_{k+1}) \\
&\qquad - A_k(\tilde{x}_{k}-\bar{\tilde{x}}_k) - B_k(U-\bar{U}) - C_k(S-\bar{S})\nonumber
\end{align} 
Note that arbitrarily chosen $\bar{\tilde{X}}$, $\bar{U}$, and $\bar{S}$ are, in general, not continuous-time feasible, i.e., $\bar{\tilde{x}}^k(\tau_{k+1}) \ne \bar{\tilde{x}}_{k+1}$, for $k=1,\ldots,N-1$. We refer to \cite[Sec. II]{malyuta2019discretization} and \cite[Sec 2.3]{kamath2023seco} for further insights and related methods for discretization and linearization of nonlinear dynamics.
%
\subsection{Pointwise Constraint Violation Bound}\label{subsec:pointwise-bound}
Given a solution for \eqref{ocp-reform-ctrlprm}, we can obtain the state trajectory $x$, control input $u$, and final time $\tfinal$. In addition, we can construct a time grid $t_1<\ldots<t_N$, satisfying $t_1 = \tinit$ and $t_N = \tfinal$. Node $t_k$ is the time instant corresponding to $\tau_k$, for $k=1,\ldots,N$. The length of interval $[t_k,t_{k+1}]$ is denoted by $\Delta t_k$, for $k=1,\ldots,N-1$, and $\Delta t_{\min}$ denotes a lower bound for the lengths of the $N-1$ intervals. 
Then we can establish an upper bound for the pointwise violation of the path constraints when \eqref{ocp-reform-ctrlprm:relax} is satisfied for a specified $\epsilon>0$, and with the choice of exterior penalties \eqref{choose-ext-penal}. 
\begin{thm} For any $i=1,\ldots,\Nineq$, $j=1,\ldots,\Neq$, and $k=1,\ldots,N-1$, there exist $\omega_{g_i}$ and $\omega_{h_j}$ such that 
\begin{align}
    \square(t,x(t),u(t)) \le{} & \delta_{\square}(\epsilon) = (4\epsilon\omega_{\square})^{\frac{1}{3}}\label{gihj-bnd}
\end{align}
whenever $\epsilon$ in \eqref{ocp-reform-ctrlprm:relax} satisfies
\begin{align}
    \epsilon\le\omega^2_{\square}\frac{\Delta t_{\min}^3}{4}
\end{align}
for $t\in[t_k,t_{k+1}]$ and $\square = g_i,h_j$.\label{thm:pointwise-bound}
\end{thm}
\begin{pf} See Appendix \ref{app:pointwise-bound}.
\end{pf}
While the analysis in Appendix \ref{app:pointwise-bound} assumes \eqref{choose-ext-penal}, other valid exterior penalty functions can be handled in a similar manner. Further, note that the bound in \eqref{gihj-bnd} is consistent because $\delta_\square(\epsilon)$ is strictly monotonic for $\epsilon\in\bRp$ and $\delta_\square(\epsilon) \to 0$ as $\epsilon \to 0$, where $\square=g_i,h_j$. We can use these bounds to select an $\epsilon$ that is numerically significant yet physically insignificant for the underlying dynamical system and constraints. 
\begin{rem}\label{rem:cnstr-tighten}
The pointwise bound for the inequality path constraints specifies the amount of constraint tightening required for the solutions obtained with relaxation \eqref{ocp-reform-ctrlprm:relax} to satisfy the constraints with no pointwise violation. More precisely,
\begin{align*}
 & \int_{t_{k}}^{t_{k+1}}\normplus{g_i(t,x(t),u(t))+\delta_{g_i}(\epsilon)}^2\mathrm{d}t \le \epsilon\\
 \implies &~g_i(t,x(t),u(t)) \le 0 \quad \forall\,t\in[t_k,t_{k+1}]
\end{align*}
We cannot, however, specify such a tightening for the equality constraints that ensures pointwise satisfaction. We can only control the extent of pointwise violation using \eqref{hj-bnd}. Specialized techniques such as the one proposed in \cite{bordalba2022direct} could be explored as a future direction.
\end{rem}
\begin{rem} \label{rem:cnstr-refine}
The estimate of $\omega_{g_i}$ and $\omega_{h_j}$ can be conservative whenever it is challenging to tightly bound $f(t,x(t),u(t))$. As a result, the bound in \eqref{gihj-bnd} can become conservative. In such cases, a solution to \eqref{ocp-reform-ctrlprm} can be iteratively refined until the pointwise violations are within a desired tolerance. We can successively solve \eqref{ocp-reform-ctrlprm} by reducing $\epsilon$ by a constant factor and warm-starting with the previous solution.  This refinement process will terminate since $\delta_{\square}$ is strictly monotonic and consistent (i.e., $\lim_{\epsilon\to0}\delta_{\square}(\epsilon) = 0$).
\end{rem}
%
The key distinction between the approach in Remark \ref{rem:cnstr-refine} and the mesh refinement techniques for alleviating inter-sample constraint violation \cite{fu2015local,fontes2019guaranteed,fu2021dynamic} is that the former does not add more node points to the discretization grid, i.e., $N$ remains fixed. 
The proposed framework can achieve continuous-time constraint satisfaction (up to a tolerance $\epsilon$) on a sparse discretization grid. Most direct methods verify the extent of inter-sample constraint violation a posteriori, i.e., after solving the trajectory optimization problem for a chosen discretization grid. In contrast, the proposed framework specifies the extent of allowable inter-sample constraint violation within the optimization problem. 
%
\section{Prox-Linear Method for Numerical Optimization}
We propose a solution method which applies an SCP algorithm called the prox-linear method \cite{lewis2016proximal,drusvyatskiy2018error} to an unconstrained minimization obtained from transforming \eqref{ocp-reform-ctrlprm} via $\ell_1$ exact penalization of \eqref{ocp-reform-ctrlprm:dyn} and \eqref{ocp-reform-ctrlprm:bc} \cite[Eq. 17.22]{nocedal2006numerical}. 
%
\subsection{Exact Penalization} \label{subsec:exact-penalty}
We can construct an unconstrained minimization problem from the constrained problem \eqref{ocp-reform-ctrlprm} by penalizing the constraints. The penalty terms have weights as tunable hyperparameters in the unconstrained problem. The penalty functions are said to be \textit{exact} whenever, for a finite penalty weight, the stationary points and local minimizers of the unconstrained problem can capture the KKT points and local minimizers of the constrained problem. A widely-used example of such penalty functions, which are typically nonsmooth, are the $\ell_1$ penalty functions \cite[Eq. 17.22]{nocedal2006numerical}. Quadratic penalty functions, in contrast, are not exact, i.e., the penalty weight must tend to infinity in order to obtain the first-order optimal points of the constrained problem. Clearly, adopting an exact penalization of constraints is practical for designing a numerical solution method since the solution to a single unconstrained problem can deliver a solution to the original constrained problem without needing to solve a sequence of problems where penalty weight is arbitrarily increased. 

Standard results on exact penalization \cite{han1979exact,di1988exactness,di1989exact}, \cite[Chap. 17]{nocedal2006numerical} for general constrained optimization appeared close to three decades ago, and are widely used. However, they don't directly apply for the proposed formulation since the construction of \eqref{ocp-reform-ctrlprm} is atypical, in that it possesses convex constraints: \eqref{ocp-reform-ctrlprm:relax}, \eqref{ocp-reform-ctrlprm:ctrlcvx} (represented with convex sets), and nonconvex constraints: \eqref{ocp-reform-ctrlprm:dyn}, \eqref{ocp-reform-ctrlprm:bc} (represented with functions). While transforming \eqref{ocp-reform-ctrlprm} into an unconstrained minimization, we wish to subject only the nonconvex constraints to exact penalization, while using indicator functions to represent the convex constraints, so that they can be directly parsed without approximations to a conic optimization solver iteratively called within the prox-linear method. The standard exact penalty results cannot be directly applied to this hybrid case. To the best of authors' knowledge, an exposition on the modifications necessary to handle the hybrid case is not available in the literature. For completeness and clarity, this section presents the exact penalty results with required modifications along with the proofs. 

Given $\gamma > 0$,  function $\Theta:\bR^{n_z}\to\bR$ is defined as 
\begin{eqboxb}
\begin{align}
\Theta_\gamma(z) ={} & \tilde{L}(\tilde{x}_N) + \indic_{\mc{Z}}(z)  \label{penalty-fun}\\
&  + \gamma\ones{}^\top\normplus{\tilde{P}(\tilde{x}_1,\tilde{x}_N)} + \gamma\|\tilde{Q}(\tinit,\tilde{x}_1,\tilde{x}_N)\|_1 \nonumber\\ 
& + \gamma\sum_{k=1}^{N-1}\|F_k(\tilde{x}_{k+1},\tilde{x}_k,U,S)\|_1  \nonumber
\end{align}
\end{eqboxb}
where $n_z = \ntildx N+n_uN_u+N_s$, $z = (\tilde{X},U,S)$, $\tilde{X} = (\tilde{x}_1,\ldots,\tilde{x}_N)$, and $\mc{Z} = \tilde{\mc{X}}\times\mc{U}\times\mc{S}$ is a closed convex set with 
$$
    \tilde{\mc{X}} = \left\{\begin{array}{l}
     (\tilde{x}_1,\ldots,\tilde{x}_N)\\
     \in\bR^{n_{\tilde{x}}N}
    \end{array}\middle|\,\begin{array}{l}
    \selector{y}(\tilde{x}_{k+1}-\tilde{x}_k) \le \epsilon,\\ \tilde{x}_k\in\bR^{\ntildx},~k=1,\ldots,N-1\end{array}\right\}
$$
The $\ell_1$-penalization of constraints \eqref{ocp-reform-ctrlprm:dyn} and \eqref{ocp-reform-ctrlprm:bc} involves penalty functions  $\square\mapsto\normplus{\square}$ and $\square\mapsto\|\square\|_1$, respectively. 

Our goal then is to compute a local minimizer of $\Theta_\gamma$ where the penalty terms are zero. We start by seeking a stationary point of $\Theta_\gamma$, since local minimizers are stationary points. Determining stationary points is a relatively simpler task than directly attempting to solve the nonconvex problem \eqref{ocp-reform-ctrlprm}. We can rewrite the indicator function for $\mc{Z}$ as explicit convex constraints\footnote{Commonly encountered convex sets such as polytope, ball, box, second-order cone, halfspace, hyperplane etc. can be directly represented in conic optimization solvers, such as ECOS \cite{domahidi2013ecos} and MOSEK \cite{mosek}, as the intersections of canonical cones: zero cone, nonnegative orthant, second-order cone, and the cone of positive semidefinite matrices.} while computing a stationary point of $\Theta_\gamma$ via SCP. A stationary point $z^\star$ satisfies $0\in\partial\Theta_\gamma(z^\star)$, i.e., there exists $\mu\in\bR^{n_z}$ and $\lambda\in\bR^{n_\lambda}$ denoted by 
\begin{align}
    \lambda ={} & (\lambda^{\tilde{P}}_1,\ldots,\lambda^{\tilde{P}}_{n_P},\lambda^{\tilde{Q}}_1,\ldots,\lambda^{\tilde{Q}}_{n_Q+1},\ldots\\
                &  \,\,\,\lambda^{F_1}_{1},\ldots,\lambda^{F_{N-1}}_{n_{\tilde{x}}})\nonumber
\end{align}
with $n_\lambda = n_P+n_Q+n_{\tilde{x}}(N-1)+1$, such that
\begin{align*}
    \lambda^{\tilde{P}}_i \in{} & \partial\normplus{\tilde{P}_i(\tilde{x}^\star_1,\tilde{x}_N^\star)}, & & i = 1,\ldots,n_P\\
    \lambda^{\tilde{Q}}_i \in{} & \partial\left|\tilde{Q}_i(\tinit,\tilde{x}^\star_1,\tilde{x}_N^\star)\right|, & & i = 1,\ldots,n_Q+1\\
    \lambda^{F_k}_i \in{} & \partial\left|F_{ki}(\tilde{x}^\star_{k+1},\tilde{x}_{k}^\star,U^\star,S^\star)\right|, & & k=1,\ldots,N-1\\
    & & & i = 1,\ldots,n_{\tilde{x}}\nonumber
\end{align*}    
\begin{align}
0_{1 \times n_z} ={} & \nabla_z L\left(\tilde{x}^\star_N\right) + \mu^\top \label{stationarity}\\
&+ \gamma \sum_{i=1}^{n_P} \lambda^{\tilde{P}}_i \nabla_z \tilde{P}_i\left(\tilde{x}_1^\star,\tilde{x}_N^\star\right)\nonumber\\
&+ \gamma \sum_{i=1}^{n_Q+1} \lambda^{\tilde{Q}}_i \nabla_z \tilde{Q}_i(\tinit,\tilde{x}_1^\star,\tilde{x}_N^\star) \nonumber\\
&+ \gamma \sum_{k=1}^{N-1} \sum_{i=1}^{n_{\tilde{x}}} \lambda^{F_k}_i \nabla_z F_{ki}(\tilde{x}_{k+1}^\star,\tilde{x}_k^\star,U^\star,S^\star)\nonumber
\end{align}
where $z^\star = (\tilde{X}^\star,U^\star,S^\star)$ with $\tilde{X}^\star = (\tilde{x}_1^\star,\ldots,\tilde{x}_N^\star)$. The scalar-valued functions $\tilde{P}_i$, $\tilde{Q}_i$, and $F_{ki}$ are elements of $\tilde{P}$, $\tilde{Q}$, and $F_k$, respectively, indexed by $i$. The gradients of $\tilde{P}_i$, $\tilde{Q}_i$, and $F_{ki}$ are with respect to $z = (\tilde{x}_1,\ldots,\tilde{x}_N,U,S)$. The convex subdifferential of $\square \mapsto |\square|_{+}$ and $\square\mapsto|\square|$ are
\begin{subequations}
\begin{align}
\partial|w|_{+} ={} & \left\{\begin{array}{cl}
\{1\}\phantom{-} & w > 0 \\
{[0,1]}\phantom{-} & w = 0 \\
\{0\}\phantom{-} & w < 0    
\end{array}\right.\label{subdiff-normplus}\\
\partial|w|_{\phantom{+}} ={} & \left\{\begin{array}{cl}
\{1\} & w > 0 \\
{[-1,1]} & w = 0 \\
\{-1\} & w < 0
\end{array}\right.
\end{align}\label{cvx-subdiff}%
\end{subequations}
for any $w\in\bR$. We refer the reader to \cite[Sec. 3.2]{drusvyatskiy2019efficiency} and the references therein, for the definition and properties of the subdifferential of a nonconvex function.
\begin{figure}
\centering
\resizebox{\linewidth}{!}{

\begin{tikzpicture}[scale=0.83]

\filldraw[fill=blue!20,draw=blue!50,line width=1pt] (-0.25,0) ellipse (4cm and 2.5cm);
\filldraw[fill=red, fill opacity=0.2,draw=red,line width=1pt,draw opacity=0.5] (-2.5,0) ellipse (4cm and 2.5cm);
\node[align=center] at (2.6,0) {Feasible\\set of \eqref{ocp-reform-ctrlprm}};
\node[align=center] at (-5.3,0) {Stationary\\points of \\$\Theta_\gamma$};
\node[align=center] at (-1.4,-1.2) {KKT\\ points of \eqref{ocp-reform-ctrlprm}};
\filldraw[fill=violet!50,draw=violet!70,line width=1pt] (-1.4,0.6) ellipse (1.6cm and 1cm);
\node[align=center] at (-1.4,0.6) {Strict local \\ minimizers \\ of \eqref{ocp-reform-ctrlprm}};


\end{tikzpicture}

}
\caption{Exact penalization ensures that, for a large enough finite $\gamma$, all stationary points of $\Theta_\gamma$ that are feasible with respect to \eqref{ocp-reform-ctrlprm} are KKT points of \eqref{ocp-reform-ctrlprm}. Furthermore, all strict local minimizers of \eqref{ocp-reform-ctrlprm} are stationary points of $\Theta_\gamma$ for a large enough $\gamma$.}
\label{stationary-kkt}
\end{figure}
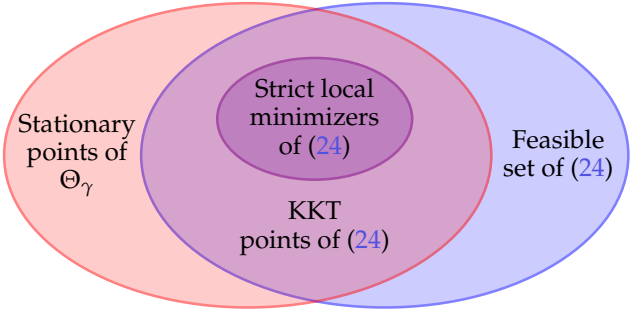

Next, we need the following lemma about the normal cone of a convex set.
\begin{lem}\label{lem:normal-cone}
Suppose that $z^\star\in\bR^n$, convex set $\mc{C} = \{z\in\bR^{n}\,|\,\check{\Omega}_i(z) \le 0,~i=1,\ldots,N_{\mc{C}}\}$ is defined with continuously differentiable convex functions $\check{\Omega}_i:\bR^n\to\bR$, and convex set $\mc{D} = \{z\in\bR^{n}\,|\,a_j^\top z + b_j = 0,~j=1,\ldots,N_{\mc{D}}\}$ is defined with vectors $a_j\in\bR^n$ and scalars $b_j\in\bR$. Assume that 1) the collection of vectors $\nabla\check{\Omega}_i(z^\star)$, for $i\in I_{\mc{C}}(z^\star)$, and $a_j$, for $j=1,\ldots,N_{\mc{D}}$ are linearly independent, where $I_{\mc{C}}(z^\star) = \{i\:|\:\check{\Omega}_i(z^\star) = 0\}$ is an index set, and 2) the intersection of $\mc{D}$ and the interior of $\mc{C}$ is nonempty. Then, we have
\begin{subequations}
\begin{align}
    \ncone_{\mc{C}}(z^\star) ={} & \Bigg\{ \sum_{i\in I_{\mc{C}}} \lambda_i\nabla\check{\Omega}_i(z^\star)^\top\,\Bigg|\,\lambda_i \in \bRp \Bigg\}\\
    \ncone_{\mc{D}}(z^\star) ={} & \Bigg\{ \sum_{j=1}^{N_{\mc{D}}}\lambda_ja_j\,\Bigg|\,\lambda_j\in\bR \Bigg\}\\
    \ncone_{\mc{C}\cap\mc{D}}(z^\star) ={} & \ncone_{\mc{C}}(z^\star) + \ncone_{\mc{D}}(z^\star)\label{sum-ncone}
\end{align}\label{ncone}%
\end{subequations}
where the set addition in \eqref{sum-ncone} is a Minkowski sum.
\end{lem}
\begin{pf}
The construction of $\ncone_{\mc{C}}(z^\star)$ and $\ncone_{\mc{D}}(z^\star)$ follow directly from \cite[Thm. 10.39, Cor. 10.44, Thm. 10.45]{clarke2013functional}, and the expression for $\ncone_{\mc{C}\cap\mc{D}}(z^\star)$ is derived from \cite[E.g. 3.5]{beck2017first} and \cite[Thm. 4.10]{clarke2013functional}, where we use the fact that $\partial\indic_{\mc{\square}}(z) = \ncone_{\square}(z)$ for any $z\in\bR^n$, where $\square=\mc{C}$, $\mc{D}$.
\end{pf}
\begin{rem}
Given a convex optimization problem with a feasible set described by $\mc{C}\cap\mc{D}$ in Lemma \ref{lem:normal-cone}, the corresponding assumptions are regularity conditions similar to LICQ and strong Slater's assumption \cite[Def. 2.3.1]{hiriart-urruty1993convex}, which are typically needed for a well-posed convex optimization problem. We require the assumptions in Lemma \ref{lem:normal-cone} on convex set $\mc{Z}$ since the SCP approach will involve a sequence of convex problems with $\mc{Z}$ as the feasible set.
\end{rem}
Then, we have the following result relating the stationary points of \eqref{penalty-fun} and KKT points of \eqref{ocp-reform-ctrlprm}.
\begin{thm}\label{thm:stationary_kkt}
Assume that $\mc{Z}$ and $z^\star$ satisfy the assumptions of Lemma \ref{lem:normal-cone}. If $z^\star$ is a stationary point of $\Theta_\gamma$ and is feasible with respect to \eqref{ocp-reform-ctrlprm}, then $z^\star$ is a KKT point of \eqref{ocp-reform-ctrlprm}. Conversely, if $z^\star$ is a KKT point of \eqref{ocp-reform-ctrlprm} and $\gamma$ is sufficiently large, then $z^\star$ is also a stationary point of $\Theta_\gamma$.
\end{thm}
%
\begin{pf}
Suppose that $z^\star$ is feasible with respect to \eqref{ocp-reform-ctrlprm}, where $z^\star = (\tilde{X}^\star,U^\star,S^\star)$ with $\tilde{X}^\star = (\tilde{x}_1^\star,\ldots,\tilde{x}_N^\star)$. Note that $z^\star$ is a KKT point of \eqref{ocp-reform-ctrlprm} if and only if there exists $\mu\in\ncone_{\mc{Z}}(z^\star)$ (using Lemma \ref{lem:normal-cone}) and $\lambda = (\lambda_1, \ldots, \lambda_{n_\lambda})\in \bR^{n_\lambda}$ such that 
\begin{align}
& \nabla_z L\left(\tilde{x}^\star_N\right) + \mu^\top + \lambda^{\top}\begin{bmatrix}\nabla_z\tilde{P}(\tilde{x}_1^\star,\tilde{x}_N^\star)\\\nabla_z\tilde{Q}(\tinit,\tilde{x}_1^\star,\tilde{x}_N^\star)\\\nabla_zF_1(\tilde{x}^\star_2,\tilde{x}^\star_1,U^\star,S^\star)\\\vdots\\\nabla_zF_{N-1}(\tilde{x}^\star_N,\tilde{x}^\star_{N-1},U^\star,S^\star)\end{bmatrix}\nonumber\\
& {}= 0_{1\times n_z} \label{kkt-stat}
\end{align}
and 
\begin{equation}
\lambda_i \tilde{P}_i(\tilde{x}_1^\star,\tilde{x}_N^\star) = 0\label{kkt-complement}
\end{equation}
for $i = 1, \ldots, n_P$. Similar to \eqref{stationarity}, the gradients of $\tilde{P}$, $\tilde{Q}$, and $F_k$ in \eqref{kkt-stat} are with respect to $z = (\tilde{x}_1,\ldots,\tilde{x}_N,U,S)$.

If $z^\star$ is a stationary point of $\Theta_\gamma$, then there exists $\mu\in\ncone_{\mc{Z}}(z^\star)$ and $\hat{\lambda} = (\hat{\lambda}_1,\ldots,\hat{\lambda}_{n_\lambda})\in\bR^{n_\lambda}$ such that \eqref{kkt-stat} holds, where
\begin{align}
\frac{1}{\gamma}\hat{\lambda}  \in{} & \partial\normplus{\tilde{P}_1} \times \ldots \times  \partial\normplus{\tilde{P}_{n_P}} \times \partial |\tilde{Q}_1| \times \ldots \label{lambda-in-subdiff}\\
 & {} \times \partial |\tilde{Q}_{n_Q+1}| \times \partial |F_{11}| \times \ldots \times \partial |F_{(N-1)n_{\tilde{x}}}|\nonumber
\end{align}
The arguments of $\tilde{P}_i$, $\tilde{Q}_i$, and $F_{ki}$, which are elements of $z^\star$, are omitted in \eqref{lambda-in-subdiff} for brevity. Observe that $\hat{\lambda}$ satisfies \eqref{kkt-complement} due to the feasibility of $z^\star$ for \eqref{ocp-reform-ctrlprm} and due to the definition of $\partial\normplus{\square}$ in \eqref{cvx-subdiff}. Therefore, $z^\star$ is a KKT point of \eqref{ocp-reform-ctrlprm}.

Conversely, if $z^\star$ is a KKT point of \eqref{ocp-reform-ctrlprm}, then there exists $\mu\in\ncone_{\mc{Z}}(z^\star)$ and $\lambda = (\lambda_1,\ldots,\lambda_{n_\lambda})$ such that \eqref{kkt-stat} and \eqref{kkt-complement} hold. Note that $z^\star$ is feasible with respect to \eqref{ocp-reform-ctrlprm} and suppose that $\gamma > \|\lambda\|_\infty$. Then, we obtain $\lambda_{n_P+i}/\gamma\in\partial |\tilde{Q}_i(\tinit,\tilde{x}_1^\star,\tilde{x}_N^\star)| = [-1,1]$, for $i = 1, \ldots, n_Q+1$, and $\lambda_{n_P+n_Q+n_{\tilde{x}}(k-1)+i+2}/\gamma\in\partial |F_{ki}(\tilde{x}_{k+1}^\star,\tilde{x}_{k}^\star,U^\star,S^\star)| = [-1,1]$, for $k=1,\ldots,N-1$, $i = 1, \ldots, n_{\tilde{x}}$.
Furthermore, for $i=1,\ldots,n_P$, we have the following. If $\tilde{P}_i(\tilde{x}_1^\star,\tilde{x}_N^\star) < 0$, then $\lambda_i/\gamma \in \partial \normplus{\tilde{P}_i(\tilde{x}^\star_1,\tilde{x}_N^\star)} = \{0\}$. Alternatively, if $\tilde{P}_i(\tilde{x}_1^\star,\tilde{x}_N^\star) = 0$, then $0 \le \lambda_i\le \gamma$, i.e., $\lambda_i/\gamma\in\partial\normplus{\tilde{P}_i(\tilde{x}^\star_1,\tilde{x}_N^\star)} = [0,1]$.

Therefore, due to \eqref{kkt-stat} and \eqref{lambda-in-subdiff}, $z^\star$ is a stationary point of $\Theta_\gamma$ for a large enough $\gamma$.%
\end{pf}
Next we relate the local minimizers of \eqref{ocp-reform-ctrlprm} to the stationary points and local minimizers of $\Theta_\gamma$.
\begin{thm}\label{thm:local-min}
Suppose LICQ holds at $z^\star$ and it is a strict local minimizer of \eqref{ocp-reform-ctrlprm}. Then, there exists $\bar{\gamma} > 0$ such that for all $\gamma \ge \bar{\gamma}$, $z^\star$ is a local minimizer of $\Theta_\gamma$. %
%
%
Moreover, there exist constants $\tilde{\gamma} > 0$ and $\tilde{\sigma} > 0$  such that, for all $\gamma \ge \tilde{\gamma}$, if $z$ is a stationary point of $\Theta_\gamma$ satisfying $\|z-z^\star\| \le \tilde{\sigma}$, then $z$ is feasible with respect to \eqref{ocp-reform-ctrlprm}.
\end{thm}
\begin{pf} Consider a modification of $\Theta_\gamma$, denoted by $\tilde{\Theta}_\gamma$, where convex constraints due to $\mc{Z}$ are also $\ell_1$-penalized instead of using the indicator function. 
The first statement follows from Theorem 4.4 in \cite{han1979exact} if we first replace $\Theta_\gamma$ with $\tilde{\Theta}_\gamma$, and note that $\tilde{\Theta}_\gamma(z) \le \Theta_{\gamma}(z)$ with equality iff $z\in\mc{Z}$. The last statement follows from Propositions 3.1 and 3.3 in \cite{di1988exactness}.
\end{pf}
Stronger versions of the exact penalization results also exist. For instance, it can be shown under certain assumptions that, given a local minimizer of \eqref{ocp-reform-ctrlprm}, there exists a neighborhood around it where the stationary points of $\Theta_\gamma$ are also local minimizers of \eqref{ocp-reform-ctrlprm}. We refer the reader to \cite{di1989exact} for further details.
%
\subsection{Prox-Linear Method}
We use the prox-linear method \cite[Sec. 5]{drusvyatskiy2018error} developed for convex-composite minimization to compute stationary points of $\Theta_\gamma$. The prox-linear method is applicable when $\Theta_\gamma$ can be fit into the following template
\begin{align}
    \Theta_\gamma(z) = J(z) + H(G(z))\label{prox-lin-template}
\end{align}
where $J$ is a proper closed convex function, $H$ is an $\alpha$-Lipschitz continuous convex function, and $G$ is potentially nonconvex and continuously differentiable with a $\beta$-Lipschitz continuous gradient.

The prox-linear method finds a stationary point $z^\star$ of $\Theta_\gamma$ by iteratively minimizing a convex approximation for $\Theta_\gamma$ given by
\begin{align}
 & \hspace{-0.2cm}\Theta^\rho_{\gamma}(z;z_k) = J(z) + H\big(G(z_k) + \nabla G(z_k) (z-z_k)\big)\label{prox-lin-cvx}\\
 & \qquad\qquad~\, + \frac{1}{2\rho}\|z-z_k\|^2 \nonumber%
\end{align}%
where $z_k$ is the current iterate, with $k \ge 1$, and $\rho$ determines the proximal term weight. The nonconvex term $H(G(z))$ is convexified by linearizing $G$ at $z_k$. We denote the unique minimizer \cite[Sec. 9.1.2]{boyd2004convex} of the strongly convex function $z \mapsto \Theta^\rho_\gamma(z;z_k)$ as $z_{k+1}$. The iterative minimization of $\Theta^\rho_\gamma$ amounts to solving a sequence of convex subproblems, which is numerically very efficient and reliable. The convergence of this iterative approach is quantified using the prox-gradient mapping
$$
    \mc{G}_\rho(z) =  \frac{1}{\rho}\big(z-\underset{z^\prime}{\mathrm{argmin}}~\Theta^\rho_\gamma(z^\prime;z)\big)
$$
This is because $\mc{G}_\rho(z)=0$ iff $z$ is a stationary point of $\Theta_\gamma$ \cite[Sec. 3.3]{drusvyatskiy2019efficiency}. Observe that $\Theta_\gamma$ defined in \eqref{penalty-fun} fits the template of \eqref{prox-lin-template}: the terminal cost function $\tilde{L}$ and the exact penalties for $\tilde{P}$, $\tilde{Q}$, and $F_k$ are represented by the composition $H\circ G$, and the indicator function $\indic_{\mc{Z}}$ is represented by $J$. More precisely, for any $z=(\tilde{x}_1,\ldots,\tilde{x}_N,U,S)\in\bR^{n_z}$, we choose
\begin{subequations}
\begin{align}
J(z) ={} & \indic_{\mc{Z}}(z)\\
H(\zeta) ={} & \zeta_{\tilde{L}} + \ones{}^\top\normplus{\zeta_{\tilde{P}}} + \|(\zeta_{\tilde{Q}},\zeta_{F})\|_1\\
G(z) ={} & \big( \tilde{L}(\tilde{x}_N),\,\gamma\tilde{P}(\tilde{x}_1,\tilde{x}_N),\,\gamma\tilde{Q}(\tinit,\tilde{x}_1,\tilde{x}_N),\\
 & ~\,F_1(\tilde{x}_2,\tilde{x}_1,U,S),\ldots,F_{N-1}(\tilde{x}_N,\tilde{x}_{N-1},U,S)\big)\nonumber
\end{align}\label{connect-prox-lin-cvx-thetgam}%
\end{subequations}
where $\zeta = (\zeta_{\tilde{L}},\zeta_{\tilde{P}},\zeta_{\tilde{Q}},\zeta_{F})$ with $\zeta_{\tilde{L}}\in\bR$, $\zeta_{\tilde{P}} \in \bR^{n_P}$, $\zeta_{\tilde{Q}}\in\bR^{n_Q+1}$, and $\zeta_{F}\in\bR^{\ntildx(N-1)}$. While minimizing \eqref{prox-lin-cvx}, we translate $\indic_{\mc{Z}}$ into explicit convex conic constraints that can be handled by
convex optimization solvers such as ECOS \cite{domahidi2013ecos}, MOSEK \cite{mosek}, and PIPG \cite{yu2022extrapolated}.
%
\begin{lem} Under Assumption \ref{asm:x-bnd}, the gradients of $\tilde{L}$, $\tilde{P}$, $\tilde{Q}$, and $F_k$ are bounded on on $\mc{Z}$.\label{lem:grad-bound-prox-lin}
\end{lem}
\begin{pf}
Functions $\tilde{L}$, $\tilde{P}$, $\tilde{Q}$, and $F_k$ are continuously differentiable. The domain of $\Theta_\gamma$ is $\mc{Z}$. From Assumption \ref{asm:x-bnd} and \cite[Chap. 4, Prop. 1.4]{clarke1998nonsmooth}, state trajectories generated with a bounded control input and dilation factor (due to \eqref{ocp-reform-ctrlprm:ctrlcvx}) are bounded. Boundedness of the dilation factor also ensures that $\tfinal$ is bounded. Then, there exists a constant, $\beta$, that bounds the norm of the gradients of $\tilde{L}$, $\tilde{P}$, $\tilde{Q}$, and $F_k$ on $\mc{Z}$.
\end{pf}
\begin{figure}[!htbp]
    \centering
    \resizebox{\linewidth}{!}{


\begin{tikzpicture}[scale=1, transform shape, node distance=1.3cm, line cap=round,
block/.style={draw, rectangle, rounded corners, align=center, minimum width=2.5cm, minimum height=1.5cm, fill=white, line width=0.875pt}, 
decision/.style={draw, diamond, rounded corners, align=center, minimum width=2.375cm, minimum height=2.375cm, fill=white, line width=0.875pt},
arrow/.style={->, line width=0.875pt}, 
circleblock/.style={draw, circle, align=center, minimum size=0.3cm, line width=0.875pt},
]


\filldraw[fill=beige!30, draw=black!30, line width=0.875pt, rounded corners, style=dashed] (11.4,-1.1) -- (11.4,-6.83) -- (-1.85,-6.83) -- (-1.85,-3.55) -- (2.75,-3.55) -- (2.75,-1.1) -- cycle;

\node[block, fill=red!20] (ocp) {Optimal Control \\ Problem};
\node[block, fill=red!20, right=of ocp] (crt) {Constraint Reformulation \\ \& Time Dilation};
\node[block, fill=red!20, right=of crt] (paramdisc) {Parameterization, \\ Discretization, \\ \& $\ell_1$ Penalization};

\node[block, fill=blue!20, below=0.875cm of ocp] (init) {Initialization};
\node[decision, fill=blue!20, right=of init, below=0.44cm of crt, text width=0.75cm, text height=0.75cm] (firstiter) {};
\node[align=center, yshift=14.5pt] at (firstiter) {{$1^{\mathrm{st}}$} \\ \small{Iteration?}};
\node[block, fill=blue!20, right=of firstiter, below=0.875cm of paramdisc] (cvx) {Convexification};

\node[decision, fill=blue!20, below=0.73125cm of firstiter, text width=0.75cm, text height=0.75cm] (conv) {};
\node[align=center] at (conv) {Converged?};

\node[block, fill=blue!20, below=1.6cm of cvx] (solve) {Solve Convex \\ Subproblem};
\node[block, fill=blue!20, below=1.6cm of init] (update) {Update \\ Hyperparameters};
\node[circleblock, fill=green!20, line width=1.25pt, below=1cm of conv, label=below:Solution] (solution) {};

\draw[arrow] (ocp) -- (crt);
\draw[arrow] (crt) -- (paramdisc);
\draw[arrow] (paramdisc) -- (cvx);

\draw[arrow] (cvx) -- (solve);
\draw[arrow, shorten >=-1.75pt] (solve) -- (solve -| conv.east);
\draw[arrow, shorten <=-1.75pt] (conv) -- node[right, xshift=-6.25pt, yshift=10pt] {No} (conv -| update.east);
\draw[arrow, shorten <=-2.5pt] (conv) ++(0,-1.23) -- node[below, xshift=-12pt, yshift=7.5pt] {Yes} (solution);

\draw[arrow] (firstiter) ++(0.05, 0) -- (firstiter -| cvx.west);
\draw[-, line width=0.875pt, color=darkgray!75, dotted, rounded corners] (firstiter) ++(-0.001, 0) -- node[below, xshift=5pt, yshift=-2.5pt] {\small{Yes}} (firstiter);
\draw[-, line width=0.875pt, color=darkgray!75, rounded corners] (firstiter) -- ++(0, -0.001) -- node[right, xshift=1.25pt, yshift=7.75pt] {\small{No}} (firstiter);
\draw[arrow, shorten >=-20pt] (init) -- (init -| firstiter.west);
\draw[arrow, rounded corners, shorten >=-20pt] (update) - ++(0, 1.5) -| (firstiter);
\filldraw[color=darkgray!75] (firstiter) circle (2pt);

\end{tikzpicture}
}
    \vspace{-1.5em}
    \caption{The proposed framework. The \coloredblock{red!20}{black} blocks show the construction of \eqref{ocp-reform-ctrlprm} and \eqref{penalty-fun}, the \coloredblock{blue!20}{black} blocks show the components of the prox-linear method, and the \coloredblock{beige!30}{black!30} block demarcates the iterative part of the algorithm.}
    \label{fig:scp}
\end{figure}
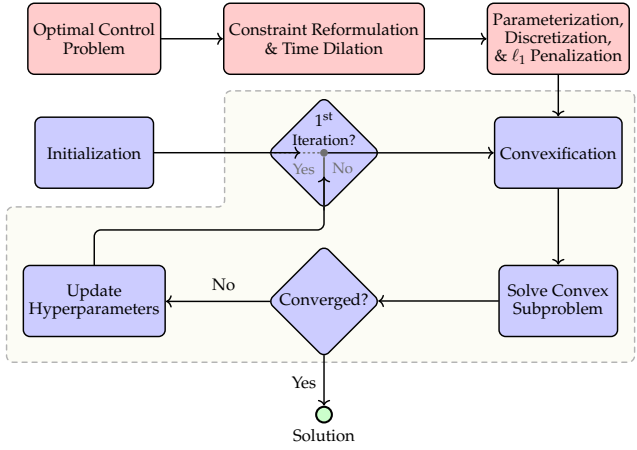
Figure \ref{fig:scp} shows the complete solution method that we propose, and Algorithm \ref{prox-linear} describes the prox-linear method. The inputs to the algorithm are the maximum number of iterations, $k_{\max}$, the termination tolerance, $\varepsilon$, and the proximal term weight, $\rho$. Note that the ``Convexification'' block in Figure \ref{fig:scp} constructs the subproblem based on \eqref{prox-lin-cvx}, which forms a convex approximation of the composition $H \circ G$. Since the convex approximation is non-affine based on our choice of $H$, we call this step ``Convexification'' rather than ``Linearization''. Furthermore, this step computes the gradient of the discretized dynamics \eqref{disc-dyn}, as described in Section \ref{subsec:jacob-dyn}.  
\begin{algorithm}[!htpb]
\caption{Prox-linear Method}
\label{prox-linear}
\begin{algorithmic}[0]
\Require $k_{\max}$, $\varepsilon$, $\rho$
\State \hspace{-0.33cm}\textbf{Initialize:} $z_1$
\State $k\gets 1$
\While{$k \le k_{\max}$ \textbf{and} $\left\|\mc{G}_\rho(z_1)\right\| > \varepsilon$}
\State $z_{k+1} \leftarrow \underset{z}{\operatorname{argmin}}~\Theta_\gamma^\rho(z;z_k)\phantom{X^{X^X}}$
\State $k \leftarrow k+1$
\EndWhile
\Ensure $z_k$
\end{algorithmic}
\end{algorithm}

We have the following result about the monotonic decrease of $\Theta_\gamma$ for a sufficiently small $\rho$.
\begin{thm}[Lemma 5.1 of \cite{drusvyatskiy2018error}]\label{thm:mono-dcrease}
At iteration $k$ of Algorithm \ref{prox-linear}, the following holds
$$
\Theta_\gamma(z_k) \geq \Theta_\gamma(z_{k+1}) + \frac{\rho}{2}(2 - \alpha \beta\rho)\|\mc{G}_\rho(z_k)\|^2
$$
Therefore, $\Theta_\gamma(z_k)$ is monotonically decreasing if $\rho \leq \frac{1}{\alpha\beta}$.
\end{thm}
The magnitude of $\mc{G}_\rho(z_k)$ is both a practical and rigorous measurement for assessing whether $z_k$ is close to a stationary point (see \cite[Thm. 5.3]{drusvyatskiy2018error} for details). We then have the following finite-stop corollary.
\begin{cor}\label{cor:near-stat}
Suppose $\Theta_{\gamma}$ is bounded below by $\Theta^\star_{\gamma}$ and $\rho \leq \frac{1}{\alpha\beta}$. For a given tolerance $\varepsilon$, prox-linear method achieves
$$
    \|\mc{G}_\rho(z_k)\|^2 < \varepsilon
\quad 
\text{for}~
    k \le \frac{2\alpha\beta(\Theta_\gamma(z_1) -  \Theta^\star_{\gamma})}{\varepsilon}
$$
\end{cor}
\begin{rem} In practice,  $\beta$ can be computed by either implementing a line-search as described in \cite[Alg. 1]{drusvyatskiy2018error} or by adopting the adaptive weight update strategy in \cite[Alg. 2.2]{cartis2011evaluation}, where the equivalence between the trust-region and prox-linear methods is also discussed. The ``Update Hyperparameters'' block in Figure \ref{fig:scp} represents the use of such update techniques. The implementation of the prox-linear method can be further enhanced using the acceleration scheme described in \cite[Sec. 7]{drusvyatskiy2019efficiency}.
\end{rem}
Note that the upper bound for $\rho$ in Theorem \ref{thm:mono-dcrease} is a sufficient condition for convergence to a stationary point. In general, if the prox-linear method converges with arbitrary user-specified weights $\rho$ and $\gamma$, then the converged solution is guaranteed to be a stationary point, due to \cite[Thm. 5.3]{drusvyatskiy2018error}. 
\begin{rem} The prox-linear method is susceptible to the \textit{crawling phenomenon} exhibited by SCP algorithms \cite{reynolds2020crawling}, wherein the iterates make vanishingly small progress towards a stationary point when they are not close to one. This phenomenon can happen when $\gamma$ is chosen to be very large, which decreases the upper bound on $\rho$ for ensuring monotonic decrease of $\Theta_\gamma$. Then, the lower bound $\frac{1}{2\alpha\beta}\|\mc{G}_{\frac{1}{\alpha\beta}}(z_k)\|$ for $\Theta_\gamma(z_k)-\Theta_\gamma(z_{k+1})$ becomes very small even if $\mc{G}_{\frac{1}{\alpha\beta}}(z_k)$ is not small. The crawling phenomenon can be remedied by scaling the decision variables and the constraint functions so that their orders of magnitude are similar.
\end{rem}
\begin{rem} The penalized trust region (PTR) SCP algorithm in \cite{szmuk2020successive,reynolds2020dual} solves a sequence convex subproblems similar to the ones generated by the prox-linear method. The PTR algorithm also implements a class of time-dilation, discretization, and parameterization operations to the optimal control problem \eqref{ocp}, albeit in an order that is different from what we propose. The resulting sequence of convex subproblems that PTR solves can fit within the template of \eqref{prox-lin-template} after minor modifications. The weights of the proximal and exact penalty terms are heuristically selected to obtain good convergence behavior.
\end{rem}
In summary, the proposed framework based on the prox-linear method globally converges to a stationary point of $\Theta_\gamma$ (Theorem \ref{thm:mono-dcrease}), and if the stationary point is feasible with respect to \eqref{ocp-reform-ctrlprm}, then it is a KKT point (Theorem \ref{thm:stationary_kkt}). Furthermore, each strict local minimizer of \eqref{ocp-reform-ctrlprm} with LICQ satisfied is a local minimizer of $\Theta_\gamma$ for a large-enough finite $\gamma$ (Theorem \ref{thm:local-min}). 
%
\section{Exploiting Convexity}
Next we consider optimal control problem \eqref{ocp} where both initial and final times $\tinit$, $\tfinal$ are fixed, the terminal cost function $L$ is convex in the terminal state, the dynamical system \eqref{ocp:dyn} is a linear nonhomogeneous ordinary differential equation, represented by
\begin{align}
\dot{x}(t) ={} & A(t)x(t) + B(t)u(t) + w(t)\label{ct-lin-sys}
\end{align}
for $t\in\tspan$, and the constraint functions in \eqref{ocp:ineq-cnstr}, \eqref{ocp:eq-cnstr}, \eqref{ocp:bc-ineq}, and \eqref{ocp:bc-eq} are jointly convex in all arguments (except time). In particular, constraint functions $h$ and $Q$ are affine in state and control input. Furthermore, we deviate from the Mayer form by explicitly specifying a running cost term in the objective with a continuously differentiable function $Y:\bRp\times\bRnx\times\bRnu\to\bR$ that is jointly convex in the state and control input. Embedding the running cost into the dynamics, as mentioned in Section \ref{subsec:ocp}, would turn the dynamical system nonlinear and hence destroy the convexity of the optimal control problem.
The resulting optimal control problem is now in the Bolza form \cite[Sec. 3.3.2]{liberzon2011calculus}. 

To obtain a finite-dimensional optimization problem similar to \eqref{ocp-reform-ctrlprm}, we first parameterize the control input and discretize \eqref{ct-lin-sys}. Time-dilation is not necessary in the fixed-final-time setting, and constraint reformulation is postponed to the end with a minor modification to exploit convexity of the constraints. We apply parameterization \eqref{aug-ctrl-param} directly to the control input in the time domain, i.e., $\nu(t) = (\Gamma_u(t)\otimes\eye{n_u})U$, for $t\in\tspan$. To discretize \eqref{ct-lin-sys}, we form a 
grid of size $N$ in $\tspan$ denoted by $\tinit=t_1<\ldots<t_N=\tfinal$, and treat the states $x_k$ at node points $t_k$ as decision variables. Then, for each $k=1,\ldots,N-1$, we exactly discretize \eqref{ct-lin-sys} via the integral form into 
\begin{align}
x^k(t) = \Phi_x(t,t_k)x_k + \Phi_u(t,t_k)U+ \phi(t,t_k) \label{dt-lin-sys}
\end{align}
for $t\in[t_k,t_{k+1}]$, where state trajectory $x^k$ generated with control input $\nu$ and initial condition $x_k$ satisfy \eqref{ct-lin-sys} $\ae$ $t\in[t_k,t_{k+1}]$. Moreover, $t\mapsto\Phi_x(t,t_k)$, $t\mapsto\Phi_u(t,t_k)$, and $t\mapsto\phi(t,t_k)$ are computed as the solution to an initial value problem\footnote{Note that, unlike in \eqref{sensitivity-ivp}, superscript ``$k$'' is absent from $\Phi_x$, $\Phi_u$, and $\phi$ since they are independent of the state trajectory $x^k$.} similar to \eqref{sensitivity-ivp}, and we denote their terminal values as $A_k = \Phi_x(t_{k+1},t_k)$, $B_k = \Phi_u(t_{k+1},t_k)$, and $w_k = \phi(t_{k+1},t_k)$.

Next, the path constraints \eqref{ocp:ineq-cnstr} and \eqref{ocp:eq-cnstr} are reformulated as follows. For each $k = 1,\dots,N-1$, define $\tilde{g}_k:\bRp\times\bR^{n_x}\times\bR^{n_uN_u}\to\bR^{n_g}$ and $\tilde{h}_k:\bRp\times\bR^{n_x}\times \bR^{n_uN_u}\to\bR^{n_h}$ by substituting in the control parameterization and the discretized form \eqref{dt-lin-sys} into $g$ and $h$
\begin{align}
& \tilde{\square}_k(t,x_k,U) = \label{cvx-cnstr-subs}\\
& \square\big(t,\Phi_x(t,t_k)x_k \!+\! \Phi_u(t,t_k)U \! +\! \phi(t,t_k),\nu(t)\big)\nonumber
\end{align}
where $\square = g,h$. Note that each scalar-valued element of $\tilde{g}_k$ and $\tilde{h}_k$ is jointly convex in $x_k$ and $U$. Then, for each $k=1,\ldots,N-1$, we obtain
\begin{subequations}
\begin{align}
& \tilde{g}_k(t,x_k,U) \le 0,~\tilde{h}_k(t,x_k,U) = 0~\forall\,t\in[t_k,t_{k+1}]\\
& \!\!\!\!\iff \xi_k(x_k,U) \le 0 \label{cvx-ctcs}
\end{align}%
\end{subequations}
where $\xi_k : \bR^{n_x}\times \bR^{n_uN_u}\to\bRp$ is defined as
\begin{align*}
& \xi_k(x_k,U) = \int_{t_k}^{t_{k+1}}\ones{}^\top \!\normplus{\tilde{g}_k(t,x_k,U)}^2 + \ones{}^\top \tilde{h}_k(t,x_k,U)^2\mathrm{d}t
\end{align*}
Similar to \eqref{ocp-reform-ctrlprm:relax}, we need to relax \eqref{cvx-ctcs} with $\epsilon>0$ in order to facilitate exact penalization.
\begin{lem} $\xi_k$ is a convex function.
\end{lem}
\begin{pf} 
Observe that $\square \mapsto \max\{0,\square\}^2$ is a convex, non-decreasing function. Then, the composition $\ones{}^\top\normplus{\tilde{g}_k}^2$ is a convex function \cite[Sec. 3.2.4]{boyd2004convex}. Furthermore, $\ones{}^\top\tilde{h}_k^2$ is a sum of convex quadratics of its arguments since $\tilde{h}_k$ is an affine function.
\end{pf}
For each $k=1,\ldots,N-1$, the running cost over $[t_k,t_{k+1}]$ can be compactly expressed with function $Y_k : \bRnx \times \bR^{n_uN_u}\to\bR$ by substituting in the control parameterization and the discretized form \eqref{dt-lin-sys} as
\begin{align}
Y_k(x_k,U) = \int_{t_k}^{t_{k+1}}Y(t,x^k(t),\nu(t))\mathrm{d}t
\end{align}
where $x^k$ is a state trajectory for \eqref{ct-lin-sys} over $[t_k,t_{k+1}]$, generated with control input $\nu$ and initial condition $x_k$.

Finally, we obtain the convex optimal control problem subject to parameterization, discretization, and reformulation \eqref{cvx-ctcs} of path constraints with $\epsilon$-relaxation.
\begin{eqbox}
\begin{subequations}
\begin{align}
\underset{X,U}{\mathrm{minimize}}~~&~L(\tfinal,x_N) + \sum_{k=1}^{N-1}Y_k(x_k,U)\\
\mathrm{subject~to}~&~x_{k+1} = A_kx_k + B_kU + w_k \label{cvx-ocp-reform-param:dyn}\\
 &~\xi_k(x_k,U) \le \epsilon \label{cvx-ocp-reform-param:path-cnstr}\\
 &~k=1,\ldots,N-1 \nonumber \\
 &~U\in\mc{U} \label{cvx-ocp-reform-param:ctrl-cnstr}\\
 &~P(\tinit,x_1,\tfinal,x_N) \le 0\\
 &~Q(\tinit,x_1,\tfinal,x_N) = 0 \label{cvx-ocp-reform-param:eq-bc}
\end{align}\label{cvx-ocp-reform-param}%
\end{subequations}
\end{eqbox}
where $X = (x_1,\ldots,x_N)\in\bR^{n_xN}$ and the compact convex set $\mc{U}$ plays a role similar to that in Remark \ref{rem:UScnstr}.

Although $\xi_k(x_k,U) \le \epsilon$ is a convex constraint, we do not possess any structural or geometric information needed to classify it as a conic constraint of the form $\hat{A}X+\hat{B}U+\hat{c}\in\mb{K}$, where $\mb{K}$ is a Cartesian product of canonical convex cones: zero cone, nonnegative orthant, second-order cone, cone of positive semidefinite matrices etc. In fact, we only have access to first-order oracles for $Y_k$ and $\xi_k$ that provide the function value and its gradient at a query point. As a consequence, we cannot directly use existing conic optimization solvers \cite{domahidi2013ecos,garstka2021cosmo,mosek} to solve \eqref{cvx-ocp-reform-param}. 

We adopt the prox-linear method after exactly penalizing \eqref{cvx-ocp-reform-param:path-cnstr}, since it is compatible with the first-order oracles for $Y_k$ and $\xi_k$. Convexity of \eqref{ocp-reform-ctrlprm} guarantees global convergence to the set of minimizers. 

Besides incorporating the $\epsilon$-relaxation in \eqref{cvx-ocp-reform-param:path-cnstr} to facilitate exact penalization, we also assume that representation \eqref{cvx-ocp-reform-param:dyn}-\eqref{cvx-ocp-reform-param:eq-bc} for the feasible set of \eqref{cvx-ocp-reform-param} satisfies the strong Slater's assumption \cite[Def. 2.3.1]{hiriart-urruty1993convex}. This assumption is necessary and sufficient for the set of Lagrange multipliers associated with a minimizer of \eqref{cvx-ocp-reform-param} to be nonempty, compact, and convex \cite[Thm. 2.3.2]{hiriart-urruty1993convex}.
%
\subsection{Exact Penalization \& Prox-Linear Method} \label{subsec:cvx-exact-penalty}
Similar to \eqref{penalty-fun}, for a given $\gamma > 0$, consider the convex function $\Theta_\gamma : \bR^{n_xN+n_uN_u} \to \bR$ defined by 
\begin{eqboxb}
\begin{align}
\Theta_\gamma(z) ={} & L(x_N)  + \sum_{k=1}^{N-1}Y_k(x_k,U) \label{cvx-penalty-fun}\\
   & + \indic_{\mc{Z}}(z) + \gamma\sum_{k=1}^{N-1}\normplus{\xi_k(x_k,U)-\epsilon}\nonumber
\end{align}
\end{eqboxb}
where $z = (X,U)$ and $\mc{Z} = \mc{X}\times\mc{U}$ with 
\begin{align}
& \mc{X} = \left\{X \in \bR^{n_xN} \middle|\begin{array}{l}
        X = (x_1,\ldots,x_N),\,x_k\in\bR^{n_x}  \\[-0.15cm]
        x_{k+1} = A_kx_k + B_kU + w_k\\[-0.15cm]
        1\le k\le N-1\\[-0.15cm]
        P(\tinit,x_1,\tfinal,x_N) = 0\\[-0.15cm]
        Q(\tinit,x_1,\tfinal,x_N) \le 0 
    \end{array}\!\!\right\}
\end{align}
is the convex set encoding all state constraints except for \eqref{cvx-ocp-reform-param:path-cnstr}, which is subject to $\ell_{1}$ penalization. Note that when $\epsilon=0$ in \eqref{cvx-penalty-fun}, i.e., when \eqref{cvx-ctcs} is not subject to relaxation, then its penalization is not exact since leads to a quadratic penalty on the violation of path constraints.
%
\begin{lem}Under Assumption \ref{asm:x-bnd}, the gradients of $Y_k$ and $\xi_k$ are bounded on $\mc{Z}$.
\label{lem:cvx-grad-bound-prox-lin}
\end{lem}
\begin{pf}
The domain of $\Theta_\gamma$ is $\mc{Z}$. Similar to Lemma \ref{lem:grad-bound-prox-lin}, we can again invoke Assumption \ref{asm:x-bnd} and \cite[Chap. 4, Prop. 1.4]{clarke1998nonsmooth} to infer that a bounded control input parameterization (due to \eqref{cvx-ocp-reform-param:ctrl-cnstr}) generates bounded state trajectories over the time interval $\tspan$. Then, there exists a constant, $\beta$, that bounds the norm of the gradients of $Y_k$ and $\xi_k$ on $\mc{Z}$.
\end{pf}

Due to the convexity of \eqref{cvx-ocp-reform-param}, we can invoke a stronger exact penalization result than those in Section \ref{subsec:cvx-exact-penalty}.

\begin{thm}\label{thm:stat=minimizer} Suppose strong Slater's assumption holds for the constraints in \eqref{cvx-ocp-reform-param}. The set of stationary points of \eqref{cvx-penalty-fun} coincides with the set of minimizers of \eqref{cvx-ocp-reform-param} when $\gamma$ is larger than the largest magnitude Lagrange multiplier associated with the minimizers of \eqref{cvx-ocp-reform-param}. 
\end{thm}
\begin{pf} The result follows from \cite[Cor. 3.2.3, Thm. 3.2.4]{hiriart-urruty1993convex}. Note that strong Slater's assumption ensures that strong duality holds. 
\end{pf}
%
\begin{thm}
Suppose strong Slater's assumption holds for the constraints in \eqref{cvx-ocp-reform-param} and $\rho\le \frac{1}{\alpha\beta}$. Then, Algorithm \ref{prox-linear} globally converges arbitrarily close to a solution of \eqref{cvx-ocp-reform-param} in a finite number of iterations.
\end{thm} 
\begin{pf} The monotonic decrease of \eqref{cvx-penalty-fun} guaranteed by Theorem \ref{thm:mono-dcrease} when $\rho \le \frac{1}{\alpha\beta}$, together with Corollary \ref{cor:near-stat} and \cite[Thm. 5.3]{drusvyatskiy2018error}, implies that Algorithm \ref{prox-linear} will get arbitrarily close to a stationary point in a finite number of iterations. Since stationary points of \eqref{cvx-penalty-fun} are identical to the minimizers of \eqref{cvx-ocp-reform-param}, due to Theorem \ref{thm:stat=minimizer}, we have that Algorithm \ref{prox-linear} will get arbitrarily close to a minimizer of \eqref{cvx-ocp-reform-param} in a finite number of iterations. 
\end{pf}
Even though prox-linear method will globally converge to a minimizer of \eqref{cvx-ocp-reform-param}, providing a good initial point $z_1$ will dramatically improve its practical convergence behavior. A solution to the closely related convex optimization problem where all path constraints are explicitly imposed as conic constraints at finitely-many time nodes will serve as a good initialization. Computing such a solution is very efficient in practice with existing convex solvers.
%
\section{Numerical Results}\label{sec:num-results}
We provide a numerical demonstration of the proposed framework with three examples: dynamic obstacle avoidance, six-degree-of-freedom (6-DoF) rocket landing, and three-degree-of-freedom (3-DoF) rocket landing with lossless convexification. The first two examples are nonconvex problems while the third one is convex. Appendices \ref{app:dyn-obstacle-avoid} through \ref{app:3DoF-rocket-landing} describe the problem formulation (in the form of \eqref{ocp}) and the algorithm parameters for each of the examples. The results in this section are generated with code in the following repository
\begin{center}
\vspace{-0.2cm}
{\footnotesize\href{https://github.com/UW-ACL/successive-convexification}{\texttt{github.com/UW-ACL/successive-convexification}}}
\vspace{-0.2cm}
\end{center}
For all examples considered, the convex subproblems of prox-linear method are QPs, which are solved using PIPG \cite{yu2022extrapolated}. We also demonstrate the real-time capability of the proposed framework by using a C codebase generated by {\scvxgen}, a general-purpose SCP-based trajectory optimization software with customized code generation, in tandem with ECOS \cite{domahidi2013ecos}.

The proximal term weight $\rho$ is determined using a line search and the exact penalty weight $\gamma$ is determined by gradually increasing it by factors of $10$ (see the discussion in \cite[p. 78]{malyuta2022convex}). Unless otherwise specified, the initial guess $z_1$ for the prox-linear method consists of a linear interpolation between the initial and final augmented states, and a constant profile for the augmented control input. The nodes of the discretization grid in $[0,1]$ for free-final-time problems, and in $\tspan$ for fixed-final-time problems, are uniformly spaced. To ensure reliable numerical performance, all decision variables and constraint functions $g_i$ and $h_j$ are scaled \cite[Sec. IV.C.2]{reynolds2020dual} to similar orders of magnitude. We refer the reader to \cite[Sec. 4.8]{betts2010practical} for a discussion on the effects of scaling in numerical optimal control. 

We compare solutions from the proposed framework with those from an implementation where path constraints are not reformulated (i.e., they are only imposed at the discretization nodes), which we refer to as the \textit{node-only} approach. We demonstrate that the proposed framework can compute continuous-time feasible solutions on a sparse discretization grid without mesh-refinement, and that the node-only approach is susceptible to inter-sample constraint violation on the same sparse grid. Handling sparse discretization grids and eliminating mesh-refinement makes the proposed framework amenable to real-time applications. Moreover, the resulting state trajectory can activate constraints between discretization time nodes, which is a unique characteristic among other direct methods.

In the plots that follow, unless otherwise stated, we adopt the following convention. Dots denote the solution from prox-linear method: black ($\bullet$) for proposed framework and blue (\textcolor{blue!50}{$\bullet$}) for node-only approach. Lines denote the continuous-time control input and state trajectory constructed using the solution from prox-linear method: black ({\Large-}) for proposed framework and blue (\textcolor{blue!50}{\Large-}) for node-only approach. (We simulate/integrate \eqref{ocp:dyn} with the parameterized control input to obtain the state trajectory.) Red lines (\textcolor{red!50}{\Large-}) denote the constraint bounds. Inter-sample violation in the solutions from the node-only approach is highlighted with magnified inset.
%
\subsection{Dynamic Obstacle Avoidance}\label{subsec:dyn-obstacle-avoid}
\begin{figure}[!htpb]
\centering
\begin{subfigure}[b]{\linewidth}
\includegraphics[width=\linewidth]{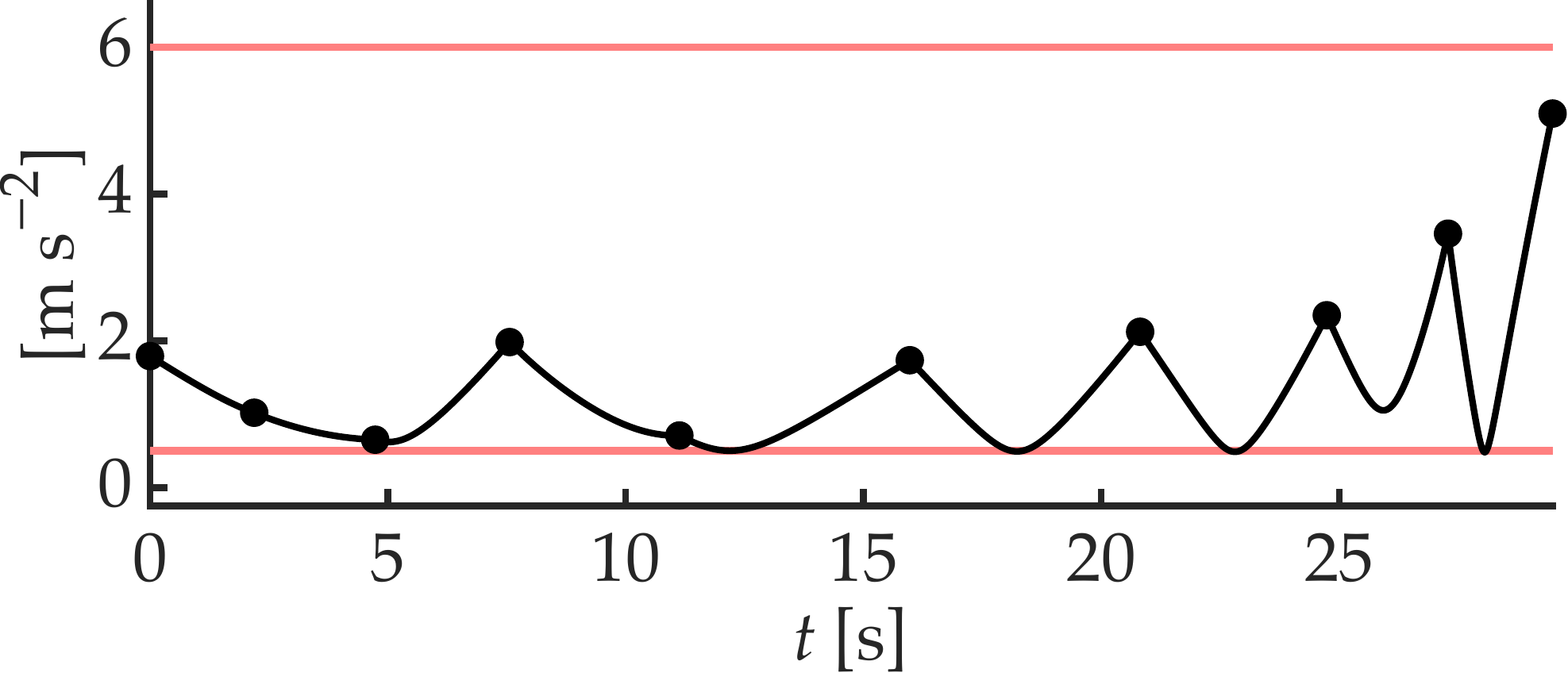}
\caption{Acceleration magnitude}
\label{fig:dyn-obs-avoid-accl}
\end{subfigure}

\vspace{0.2cm}

\begin{subfigure}[b]{\linewidth}
\includegraphics[width=\linewidth]{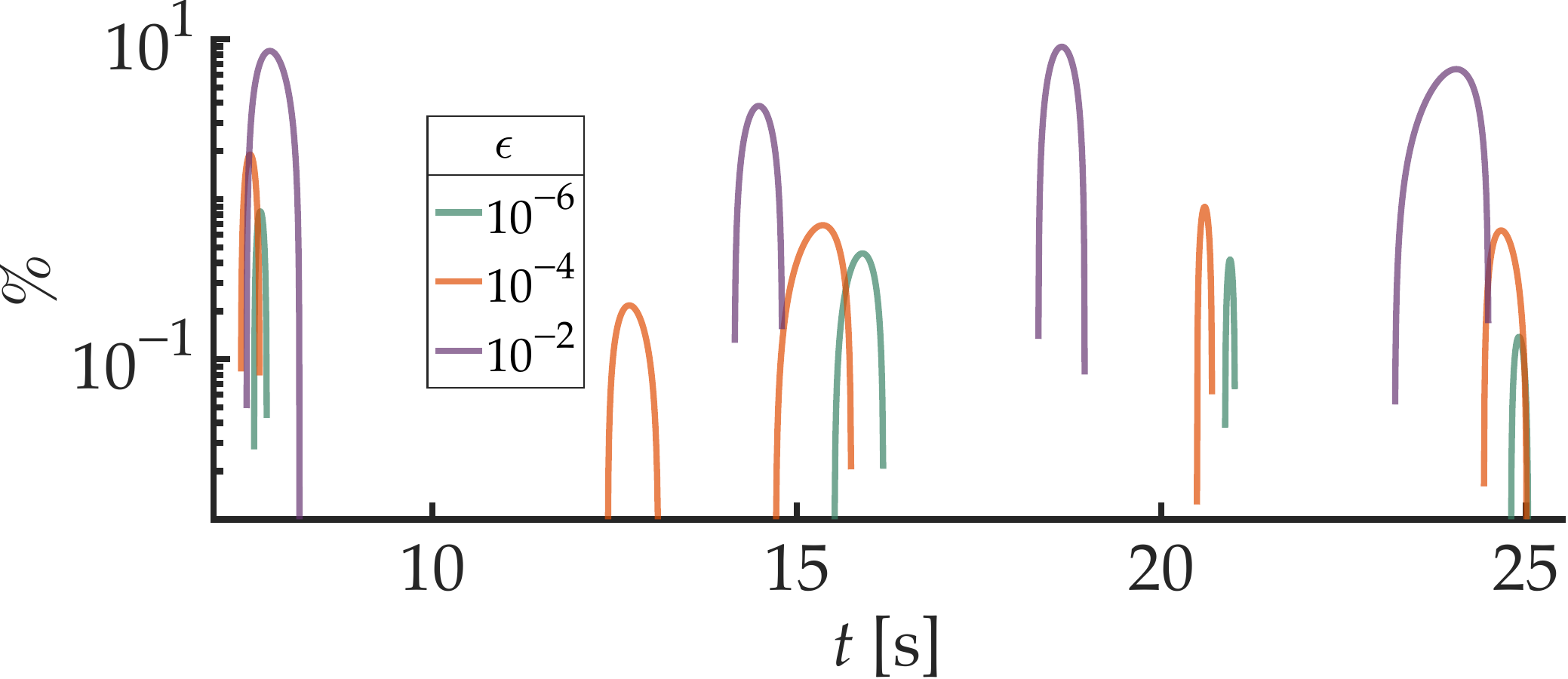}
\caption{Constraint violation}
\label{fig:dyn-obs-avoid-eps-cnstr-sweep}
\end{subfigure}
\caption{Dynamic obstacle avoidance}
\label{fig:dyn-obs-avoid}
\end{figure}
\begin{figure}[!htpb]
\centering
\begin{subfigure}[b]{\linewidth}
\includegraphics[width=\linewidth]{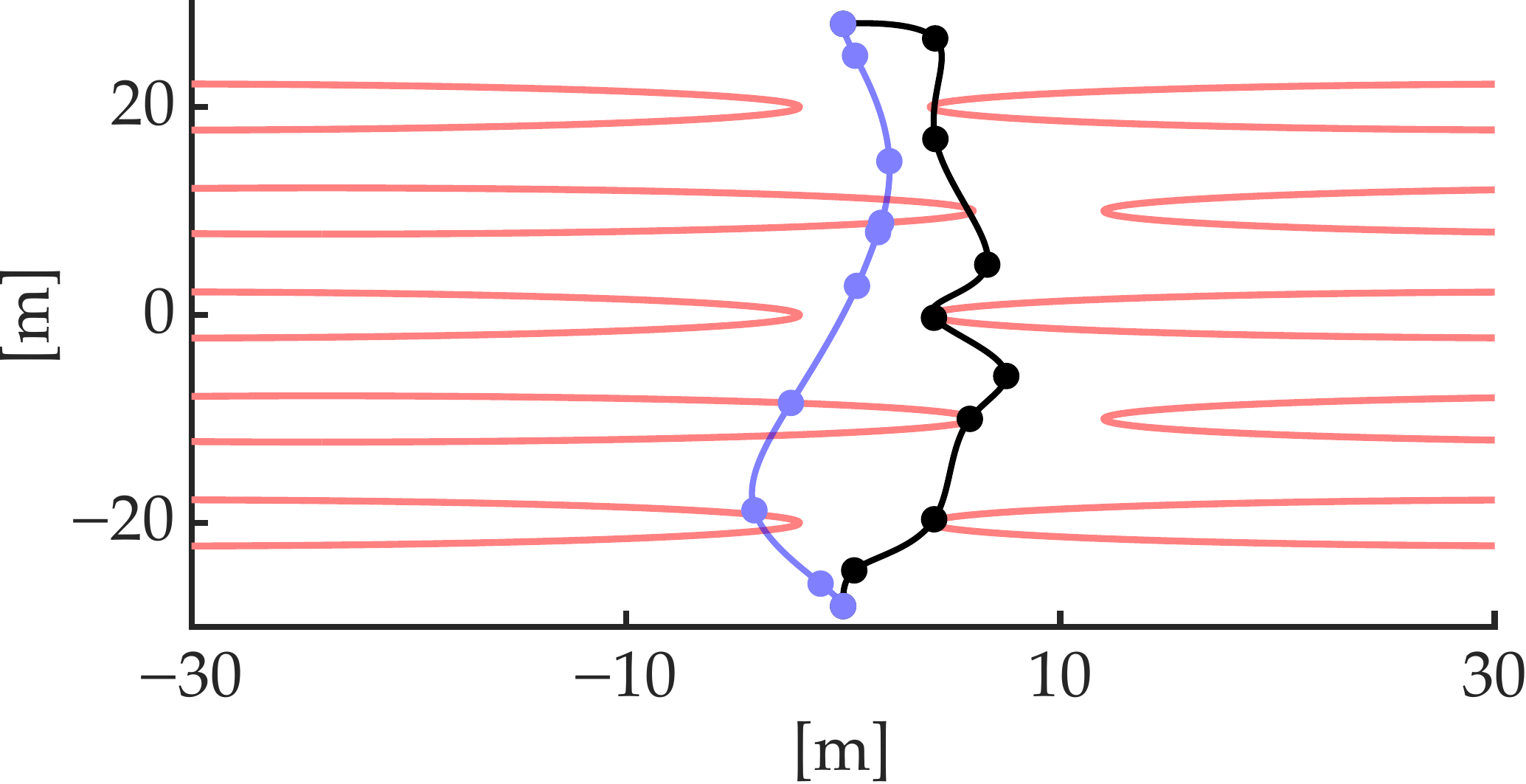}    
\caption{Position}
\label{fig:stat-obs-avoid-pos}
\end{subfigure}

\vspace{0.2cm}

\begin{subfigure}[b]{\linewidth}
\includegraphics[width=\linewidth]{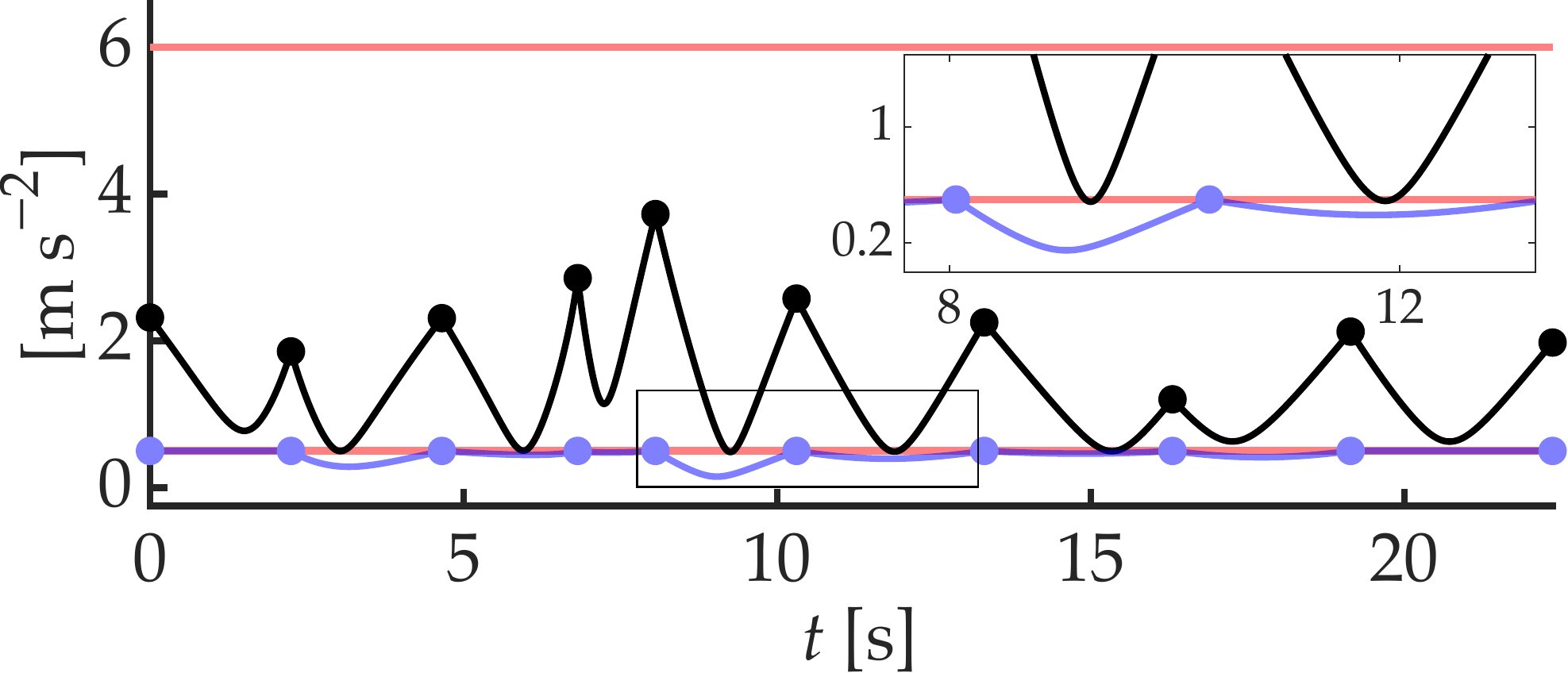}    
\caption{Acceleration magnitude}
\label{fig:stat-obs-avoid-accl}
\end{subfigure}
\caption{Static obstacle avoidance}
\label{fig:stat-obs-avoid}
\end{figure}
We consider two scenarios: one where the obstacles exhibit periodic motion and the other where they are static. Figure \ref{fig:dyn-obs-avoid-accl} shows the magnitude of control input (acceleration) in the case with dynamic obstacles obtained with the proposed framework. An animation of the corresponding position trajectory is provided in the code repository. A solution from the node-only approach is not shown in Figure \ref{fig:dyn-obs-avoid-accl}, since it fails to provide a physically meaningful solution on the same discretization grid due to aliasing in the obstacle motion. Due to the time variation in the obstacle positions, a dense discretization grid is required with the node-only approach. However, the number of distinct obstacle avoidance constraints imposed in the node-only approach grows with the discretization grid size and the number of obstacles, which is computationally expensive. On the other hand, the number of constraints imposed in the proposed framework only depends on the discretization grid size.

Figure \ref{fig:dyn-obs-avoid-eps-cnstr-sweep} shows the percentage of violation in the obstacle avoidance constraint by the simulated state trajectories obtained for different values of the constraint relaxation tolerance $\epsilon$. As long as the constraint functions are well-scaled, a state trajectory with physically insignificant continuous-time constraint violation can be obtained by picking $\epsilon$ to be several orders of magnitude greater than machine precision, which is essential for reliable numerical performance. In practice, picking $\epsilon$ close to $10^{-4}$ ensures that the continuous-time constraint violation does not exceed $1\%$. 


To compare the proposed and the node-only approaches, we turn to the case of static obstacles: Figure \ref{fig:stat-obs-avoid-pos} shows the position and Figure \ref{fig:stat-obs-avoid-accl} shows the acceleration magnitude. Note that the acceleration magnitude in the solution from the node-only approach is smaller than the lower bound (strictly smaller between the discretization nodes). As a result, the terminal state cost from the node-only approach, 5.06, is significantly smaller than that from the proposed framework, 47.91---the node-only approach optimizes the cost at the expense of inter-sample constraint violation.
%
%
\subsection{6-DoF Rocket Landing}\label{subsec:rocket-land-6DoF}
\begin{figure}[!htpb]
\centering
\begin{subfigure}[b]{\linewidth}
\includegraphics[width=\linewidth]{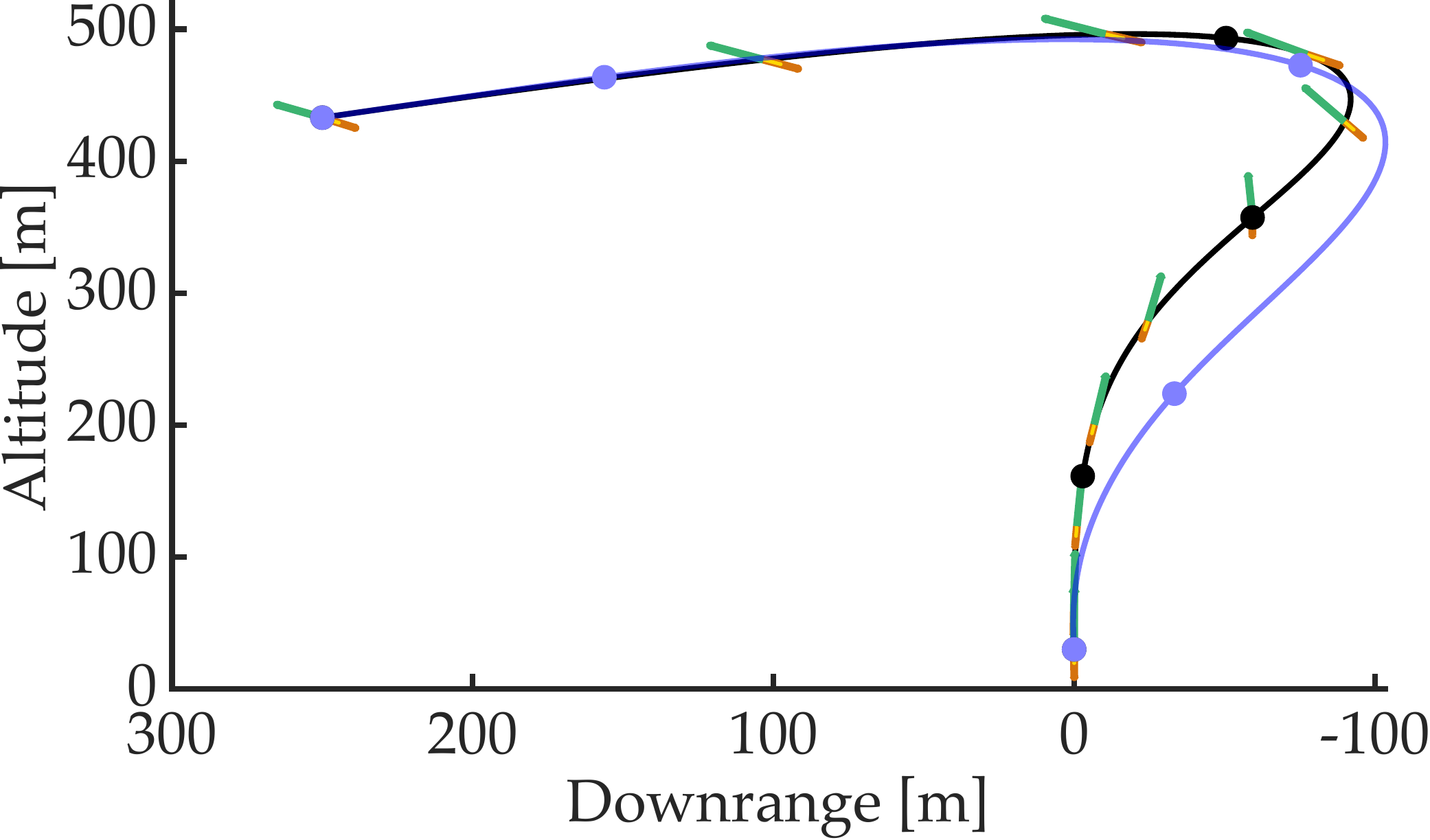}    
\caption{Position}
\label{fig:6DoF-rocket-landing-pos}
\end{subfigure}

\vspace{0.2cm}

\begin{subfigure}[b]{\linewidth}
\includegraphics[width=\linewidth]{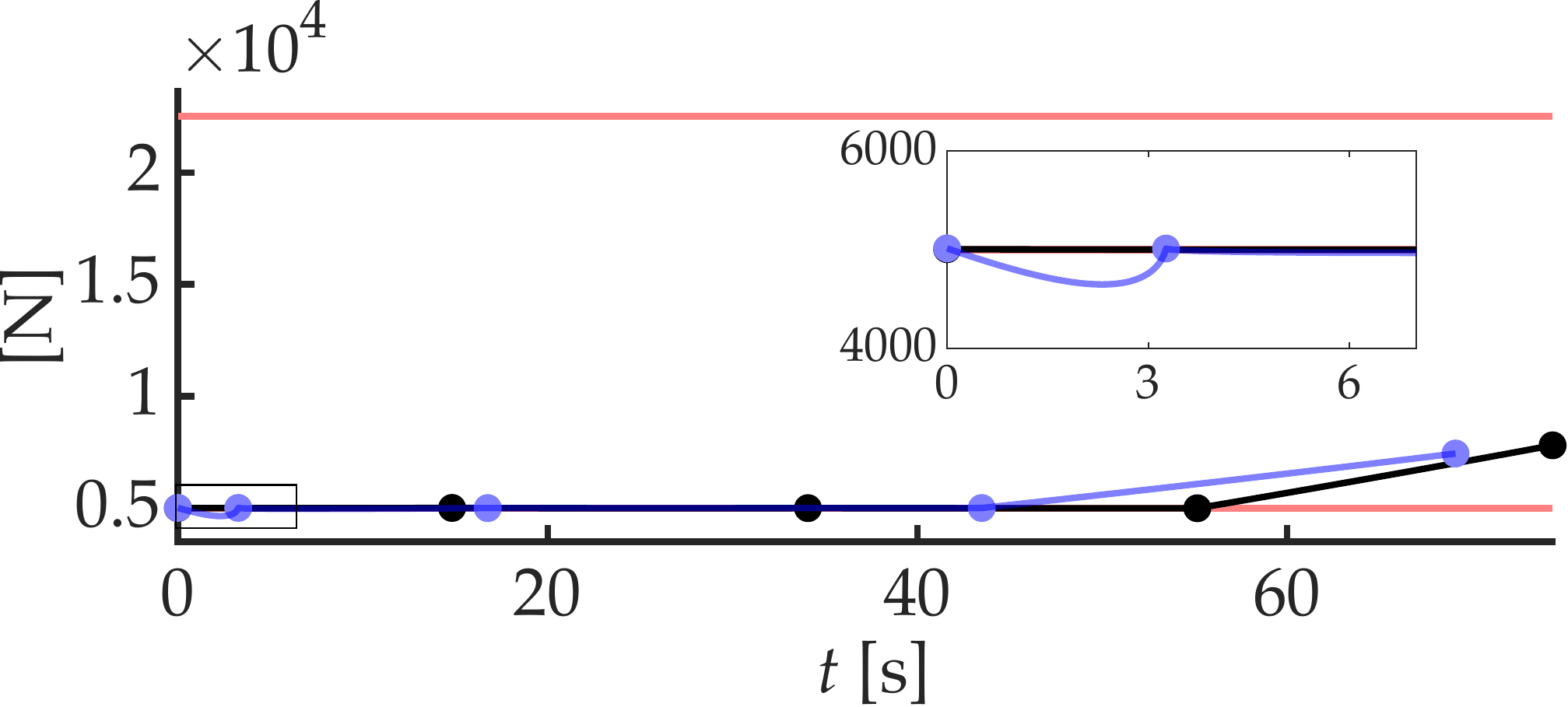}    
\caption{Thrust magnitude}
\label{fig:6DoF-rocket-landing-thrust}
\end{subfigure}

\vspace{0.2cm}

\begin{subfigure}[b]{\linewidth}
\includegraphics[width=\linewidth]{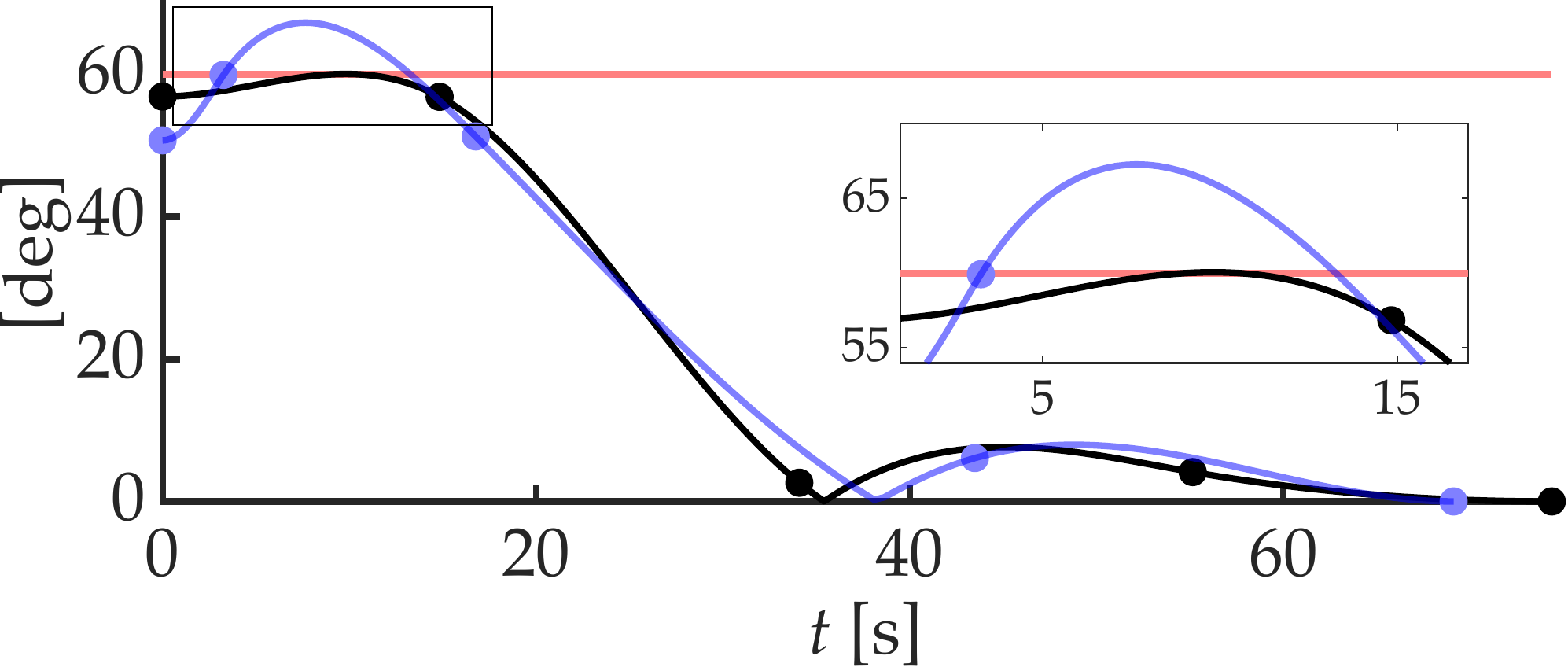}
\caption{Tilt angle}
\label{fig:6DoF-rocket-landing-tilt}
\end{subfigure}
\caption{6-DoF rocket landing}
\label{fig:6DoF-rocket-landing}
\end{figure}
Figure \ref{fig:6DoF-rocket-landing-pos} shows the position of the 6-DoF rocket along with the attitude of its body-axis (in green) and the thrust vector (as an orange-yellow plume). Figure \ref{fig:6DoF-rocket-landing-thrust} shows the thrust magnitude and Figure \ref{fig:6DoF-rocket-landing-tilt} shows the tilt of the body axis---both of which show inter-sample violation with the node-only approach.


Similar to the obstacle avoidance example, the terminal state cost with the node-only approach, $170.1$ kg, is smaller than that from the proposed framework, $180.5$ kg.
%
\subsection{3-DoF Rocket Landing (Lossless Convexification)}\label{subsec:rocket-land-3DoF}
\begin{figure}[!htpb]
\centering
\begin{subfigure}[b]{\linewidth}
\includegraphics[width=\linewidth]{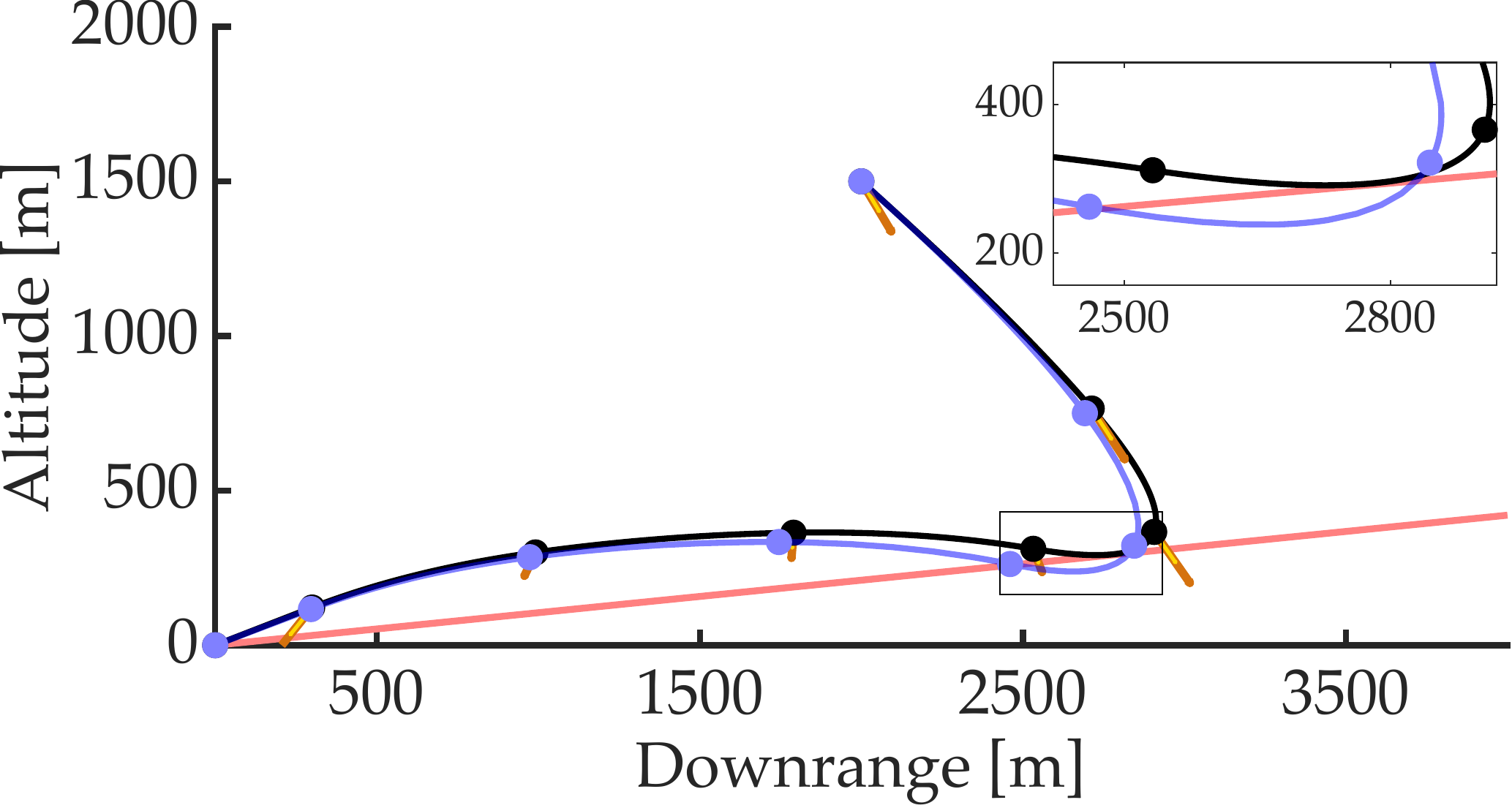}
\caption{Position}
\label{fig:3DoF-rocket-landing-pos}
\end{subfigure}

\vspace{0.2cm}

\begin{subfigure}[b]{\linewidth}
\includegraphics[width=\linewidth]{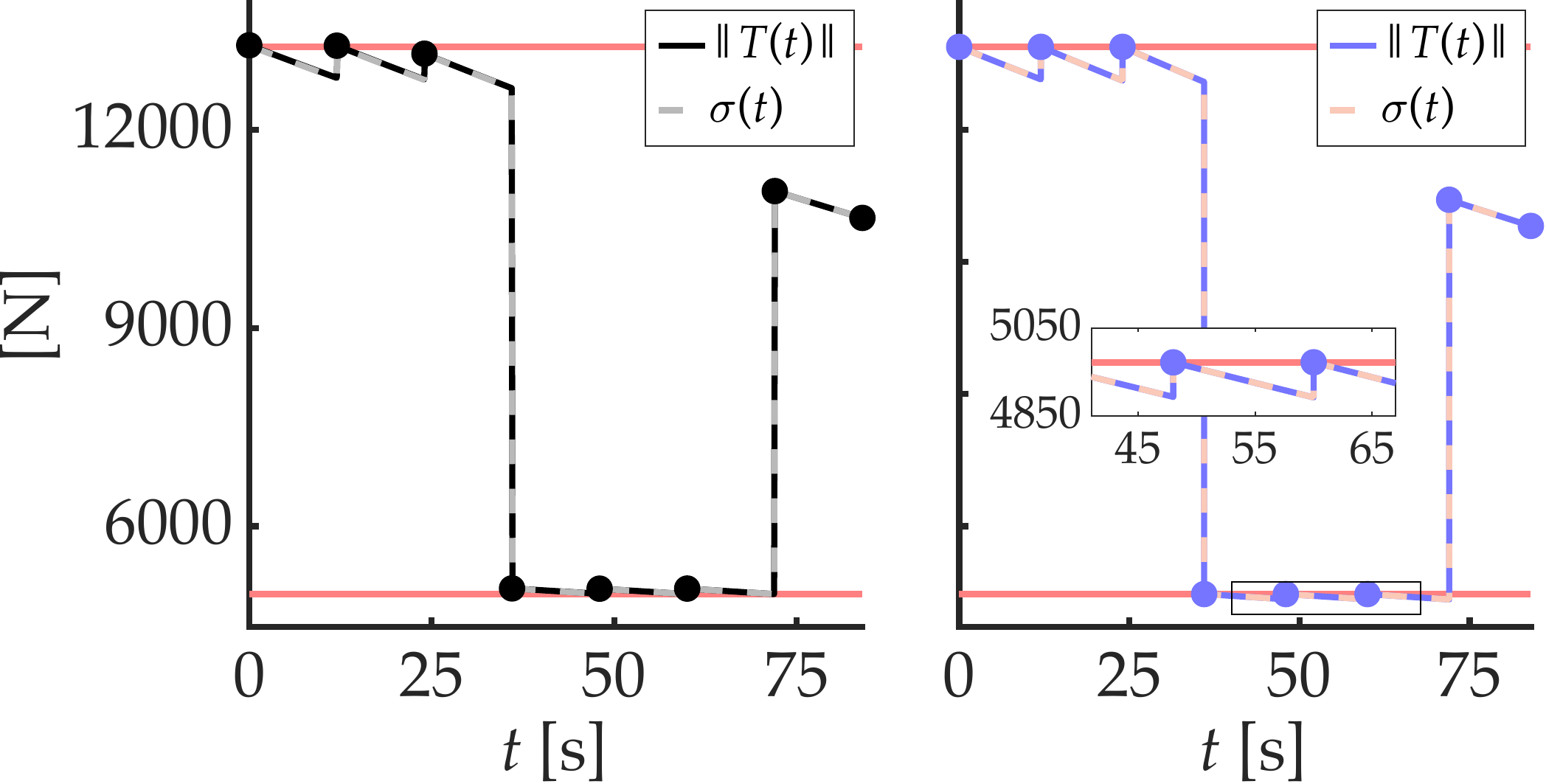}
\caption{Thrust magnitude}
\label{fig:3DoF-rocket-landing-thrust}
\end{subfigure}
\caption{3-DoF rocket landing with lossless convexification}
\label{fig:3DoF-rocket-landing}
\end{figure}
This example demonstrates the application of the proposed framework for solving a convex optimal control problem. The node-only approach, in this case, amounts to solving a single convex problem where all constraints are imposed via canonical convex cones in a convex optimization solver; we use ECOS \cite{domahidi2013ecos}. 

The solution from the node-only approach is used to warm-start the proposed framework, which then ``cleans up'' the inter-sample constraint violations. As with the previous examples, we see that the trajectory cost (in this case, fuel consumption) with the node-only approach, $351$ kg, is smaller than that from the proposed framework, $352.4$ kg. Figure \ref{fig:3DoF-rocket-landing-pos} shows the position of the point-mass rocket along with the thrust vector (as an orange-yellow plume). The glideslope constraint shows inter-sample violation with the node-only approach. Figure \ref{fig:3DoF-rocket-landing-thrust} confirms that lossless convexification holds for the solutions from both methods---proposed on the left and node-only on the right. The thrust lower bound constraint shows inter-sample violation with the node-only approach, since a zero-order-hold (ZOH) parameterization is used for the mass-normalized thrust. 

The node-only approach requires a dense discretization grid to ensure validity of the lossless convexification result \cite{acikmese2013lossless}, which was proven for the continuous-time problem (prior to discretization). The proposed framework, in contrast, provides a valid solution despite using a sparse discretization grid. Future work will examine whether the continuous-time lossless convexification result can apply directly to \eqref{cvx-ocp-reform-param}.

%
\subsection{Real-Time Performance}\label{subsec:real-time}
We provide real-time performance statistics for solving the two nonconvex examples in Table \ref{tab:real-time}. These results were obtained by executing pure C code generated using the {\scvxgen} software. The binary executables are provided in the code repository.

\begin{table}[!htpb]
\caption{\\[-0.75em]\textit{The mean solve-time and standard deviation (S.D.) over 1000 solves, along with the number of prox-linear iterations to convergence (Iters.), for the C implementation.\\[-0.25em]}}\label{tab:real-time}
{\renewcommand{\arraystretch}{1.1}
\begin{tabular}{l|l|l|l}
\hline
Problem & Solve-time & S.D. & Iters. \\\hline
Obstacle avoidance & $59.4$ ms & $1.4$ ms & 23 \\
6-DoF rocket landing & $42.5$ ms & $1.3$ ms & 14 \\
\hline
\end{tabular}}
\end{table}
We note that these solve-times are on the same order-of-magnitude as the solve-times reported in the literature for real-time powered-descent guidance \cite{dueri2017customized,reynolds2020real,elango2022customized,kamath2023seco,kamath2023customized} and real-time quadrotor trajectory optimization \cite{szmuk2017convexification,szmuk2018real,szmuk2019real,yu2023real}, and thus demonstrate that the proposed framework is amenable to online trajectory optimization for onboard/embedded applications.
%
%
\section{Conclusions}
We propose a novel trajectory optimization method that ensures continuous-time constraint satisfaction, guarantees convergence, and is real-time capable. The approach leverages an SCP algorithm called the prox-linear method along with $\ell_1$ exact penalization. When the optimal control problem is convex, we show stronger convergence property of the prox-linear method.

Future work will consider assumptions weaker than LICQ for exact penalty to hold, with the goal of developing a solution method that directly handles \eqref{disc-cnstr-bc} without relaxation, and also releasing {\scvxgen}, the in-house-developed code-generation software for general-purpose real-time trajectory optimization that was used to generate pure C code for the examples considered in the paper.
\begin{ack}
The authors gratefully acknowledge Mehran Mesbahi, Skye Mceowen, and Govind M. Chari for helpful discussions and their feedback on the initial draft of the paper.
\end{ack}

{\RaggedRight%
\bibliographystyle{unsrturl}      
\bibliography{references}%
}

\clearpage

\section*{Appendix}

\appendix

\section{Proof of Theorem \ref{thm:pointwise-bound}}\label{app:pointwise-bound}
As a consequence of the relaxation in \eqref{ocp-reform-ctrlprm:relax}, we wish to bound $\normplus{g_i}$ and $|h_j|$ on the interval $[t_k, t_{k+1}]$, for any $i=1,\ldots,\Nineq$ and $j=1,\ldots,\Neq$.
\begin{lem}
Time derivative of path constraint functions $g_i$ and $h_j$ are bounded almost everywhere on the state trajectory $x$ and control input $u$.
\end{lem}
\begin{pf}
Note that the augmented control input parameterization \eqref{aug-ctrl-param} is a piecewise polynomial, which is differentiable almost everywhere in $[0,1]$. Then, due to \eqref{ocp-reform-ctrlprm:ctrlcvx} and Remark \ref{rem:UScnstr}, $u(t)$ and $s(\tau)$ are bounded for $t\in\tspan$, $\tau\in[0,1]$, and $\dot{u}(t) = \derv{u}(t(\tau))/s(\tau)$ is well-defined and bounded $\ae$ $\tau\in[0,1]$. The state trajectory $x$ is bounded due to Assumption \ref{asm:x-bnd}. Further, $\Delta t_k$, for $k=1,\ldots,N-1$, is bounded. Consequently, $f$, $\inlpbyp{g_i}{\scriptscriptstyle\square}$, and $\inlpbyp{h_j}{\scriptscriptstyle\square}$ are bounded, where $\square = t$, $x$, and $u$. Then, through the chain-rule
\begin{align}
\dot{\square} ={} & \inlpbyp{\square}{t} + \inlpbyp{\square}{x}\dot{x} + \inlpbyp{\square}{u}\dot{u}\label{hjgi-bnd-chain-rule}
\end{align} 
with $\square=g_i,\,h_j$, we infer that $\dot{g}_i$ and $\dot{h}_j$ are bounded almost everywhere.
\end{pf}
We denote the upper bounds on the absolute values of $\dot{g}_i$ and $\dot{h}_j$ (which hold almost everywhere) with $\omega_{g_i}$ and $\omega_{h_j}$, respectively.

When \eqref{ocp-reform-ctrlprm:relax} and the penultimate row of \eqref{ocp-reform-ctrlprm:dyn} are satisfied, for any $k=1,\ldots,N-1$, we have that 
\begin{align*}
&\int_{\tau_{k}}^{\tau_{k+1}} s(\tau)\,\ones{}^\top\normplus{g(t(\tau),x(t(\tau)),u(t(\tau)))}^2\\
& \quad~~\, +s(\tau)\,\ones{}^\top h(t(\tau),x(t(\tau)),u(t(\tau)))^2\,\mathrm{d}\tau \le \epsilon\\
&\hspace{-0.3cm}\iff\!\!\!\!\int_{t_k}^{t_{k+1}}\!\!\ones{}^\top\normplus{g(t,x(t),u(t))}^2\!+\!\ones{}^\top h(t,x(t),u(t))^2\,\mathrm{d}t \le \epsilon\\
&\hspace{-0.3cm}\implies{}\!\!\!\!\int_{t_k}^{t_{k+1}}\!\normplus{g_i(t,x(t),u(t))}^2\mathrm{d}t \le \epsilon\\
&\phantom{\hspace{-0.3cm}\implies{}}\!\!\!\int_{t_k}^{t_{k+1}}\!h_j(t,x(t),u(t))^2\mathrm{d}t \le \epsilon
\end{align*}
for any $i=1,\ldots,n_g$ and $j=1,\ldots,n_h$, where we use the change of variable $t(\tau) = t_k + \int^{\tau}_{\tau_k}s(\theta)\mathrm{d}\theta$, for $\tau\in[\tau_k,\tau_{k+1}]$.
\begin{figure}[!htpb]
\begin{subfigure}[b]{0.49\columnwidth}
    \centering


\begin{tikzpicture}[scale=0.5, line cap=round]
    \pgfplotsset{
        every axis/.append style={
            xlabel={},
            ylabel={},
            xmin=-2, xmax=2,
            ymin=-20, ymax=2,        
            ylabel style={
                rotate=-90, 
                anchor=west,
                yshift=0cm,         
                xshift=-0.5cm       
            },            
            ytick=\empty,           
            grid=none,              
            axis lines=left,        
            enlargelimits=true,     
            axis line style={-To, line width=0.6pt},
            samples=400
        }
    }
    \begin{axis}[
        name=plot1,
        ytick={-0.3},
        yticklabel={$\hat{h}_j^2$},
        xtick={-2,-0.69,1.68},         
        xticklabels={$t_{k}^{\vphantom{\prime}}$, $t_{\max}^{\vphantom{\prime}}$, $t_{k+1}^{\vphantom{\prime}}$} 
        ]
        \addplot[blue, thick, line width=0.8pt, domain=-2:2, samples=100] {-1.6*(x+0.69)^2-0.4};
        \addplot[name path=curve, red, thick, line width=0.8pt, domain=-0.69:1.3, samples=100] {-11*x-8};        
        \path[name path=axis] (axis cs:-0.69,-22.2) -- (axis cs:1.3,-22.2);
        \addplot[red!20] fill between[of=curve and axis];    
        \node[above] at (axis cs:1,-3) {$h_j[t]^2$};
        \node[above] at (axis cs:-0.1,-16) {$\eta(t)$};
        \draw [black, dashed] (-0.69,-0.3) -- (axis cs:-3,-0.3);
    \end{axis}
\end{tikzpicture}

    \caption{$\hat{h}_j \le \Delta t_k \omega_{h_j}$}
    \label{fig:pointwise-a}    
\end{subfigure}
\hfill
\begin{subfigure}[b]{0.49\columnwidth}
    \centering

\begin{tikzpicture}[scale=0.5, line cap=round]
    \pgfplotsset{
        every axis/.append style={
            xlabel={},
            ylabel={},
            xmin=-2, xmax=2,
            ymin=-20, ymax=2,        
            ylabel style={
                rotate=-90, 
                anchor=west,
                yshift=0cm,         
                xshift=-0.5cm       
            },            
            ytick=\empty,           
            grid=none,              
            axis lines=left,        
            enlargelimits=true,     
            axis line style={-To, line width=0.6pt},
            samples=400
        }
    }
    \begin{axis}[
        name=plot2,
        ytick={-0.3},
        yticklabel={$\hat{h}_j^2$},
        xtick={-2,-0.75,0.65,1.68},         
        xticklabels={{$t^{\phantom{\prime}}_{k}$}, {$t^{\phantom{\prime}}_{\max}$}, {$t^\prime_{\max}$}, {$t^{\phantom{\prime}}_{k+1}$}} 
        ]
        \addplot[blue, thick, line width=0.8pt, domain=-2:2, samples=100] {-1*(x+0.75)^2-0.5};
        \addplot[name path=curve, red, line width=0.8pt,  domain=-0.75:0.65, samples=100] {-6*x-4.9};        
        \path[name path=axis] (axis cs:-0.75,-22.2) -- (axis cs:0.65,-22.2);
        \addplot[red!20] fill between[of=curve and axis];
        \node[above] at (axis cs:1,-3) {$h_j[t]^2$};
        \node[above] at (axis cs:0,-12) {$\eta(t)$};
        \draw [black, dashed] (-0.75,-0.3) -- (axis cs:-3,-0.3);
    \end{axis}    
\end{tikzpicture}

    \caption{$\hat{h}_j > \Delta t_k \omega_{h_j}$}
    \label{fig:pointwise-b}    
\end{subfigure}
\caption{Approximation of the area under $h_j[t]^2$.}
\label{fig:pointwise}
\end{figure}
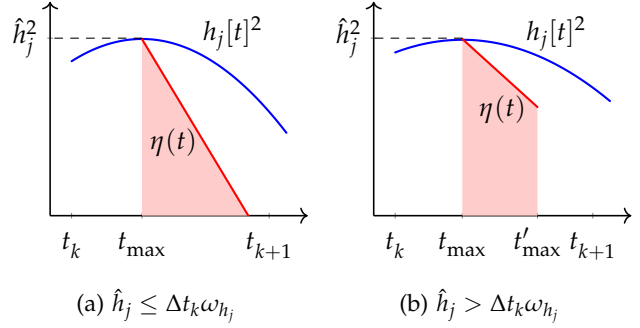

First, we examine the pointwise bound on $h_j$. We denote $h_j(t,x(t),u(t))$ compactly with $h_j[t]$, and let $\hat{h}_{j}$ denote the maximum absolute value of $h_j$ within the interval $[t_k, t_{k+1}]$. We assume, without loss of generality, that $\hat{h}_{j}$ is attained at $t_{\max} \in [t_k, t_k + 0.5\Delta t_k]$. Next, we approximate the area under $h_j[t]^2$, for $t\in[t_k,t_{k+1}]$, via an auxiliary function $\eta:[t_{\max},t_{k+1}] \to \bR$  given by
$$ \eta(t) = \normplus{\hat{h}_{j}^2 - 2 \omega_{h_j}\hat{h}_{j}(t-t_{\max})} $$ 
Note that, for $t\in[t_{\max},t_{k+1}]$, we have  
\begin{align*}
 & h_j[t]^2 - \hat{h}^2_j + 2\omega_{h_j}\hat{h}_j(t-t_{\max})\\
 \ge{} & h_j[t]^2 - \hat{h}^2_j + 2\hat{h}_j(\hat{h}_j-h_j[t])\\
   ={} & h_j[t]^2 + \hat{h}^2_j - 2\hat{h}_jh_j[t] = (h_j[t] - \hat{h}_j)^2\\
 \ge{} & 0
\end{align*}
Therefore, $h_j[t]^2 - \eta(t)$ is positive for $t > t_{\max}$ and  $h_j[t_{\max}]^2 = \eta (t_{\mathrm{max}})$. If $\hat{h}_{j} \leq \Delta t_k \omega_{h_j}$, we have that $\eta(t_{\max}^\prime) = 0$, where $t_{\max}^\prime = t_{\max} + 0.5\Delta t_k$. The area under $\eta(t)$, for $t\in [t_{\max}, t_{k+1}]$, is triangular (Figure \ref{fig:pointwise-a}). Then, 
\begin{subequations}
\begin{align}
&\hspace{-0.3cm}\epsilon\ge\int_{t_k}^{t_{k+1}} h_j[t]^2 \mathrm{d}t \ge \int_{t_{\max}}^{t_{k+1}} \eta(t)\mathrm{d}t \label{hjbnd-1}\\
&\hspace{-0.3cm}\hphantom{\epsilon}=\int_{t_{\max}}^{t_{\max} +{} \textstyle \frac{\hat{h}_{j}}{2\omega_{h_j}}} \hat{h}_{j}^2-2\omega_{h_j}\hat{h}_{j}(t-t_{\max})\mathrm{d}t\\
&\hspace{-0.3cm}\hphantom{\epsilon}=\int_{0}^{\textstyle\frac{\hat{h}_{j}}{2\omega_{h_j}}}\hat{h}_{j}^2-2\omega_{h_j}\hat{h}_{j}\theta\,\mathrm{d}\theta = \frac{\hat{h}_{j}^3}{4\omega_{h_j}}\label{pointwise-bound-case-1}
\end{align}
\end{subequations}
Conversely, if $\hat{h}_{j} > \Delta t_k \omega_{h_j}$, we have that $\eta(t_{\max}^\prime) > 0$. The area under $\eta(t)$ within $[t_{\max},t_{\max}^\prime]$ is trapezoidal (Figure \ref{fig:pointwise-b}). Then,
\begin{subequations}
\begin{align}
&\hspace{-0.1cm}\epsilon\ge\int_{t_{k}}^{t_{k+1}} h_j[t]^2\mathrm{d}t \ge \int_{t_{\max}}^{t_{k+1}} \eta(t)\mathrm{d}t\\
&\hspace{-0.1cm}\hphantom{\epsilon}>\int_{t_{\max}}^{t_{\max}+{}\textstyle\frac{\Delta t_k}{2}}\hat{h}_j^2-2\omega_{h_j}\hat{h}_j(t-t_{\max})\mathrm{d}t\\
&\hspace{-0.1cm}\hphantom{\epsilon}=\int_{0}^{\textstyle\frac{\Delta t_k}{2}}\!\!\hat{h}_{j}^2 - 2 \omega_{h_j}\hat{h}_{j} \theta\,\mathrm{d}\theta\\
&\hspace{-0.1cm}\hphantom{\epsilon}= \hat{h}^2_j\frac{\Delta t_k}{2} - \omega_{h_j}\hat{h}_j\frac{\Delta t_k^2}{4} >\hat{h}_{j}^2 \frac{\Delta t_k}{4} > \omega_{h_j}^2\frac{ \Delta t_{\min}^3}{4}
\end{align}
\end{subequations}
where $\Delta t_{\min} > 0$ is the lower bound for $\Delta t_k$; it exists because the dilation factor is positive and $\tau_{k+1}-\tau_{k}>0$, for $k=1,\ldots,N-1$. Note that $\hat{h}_{j} > \Delta t_k \omega_{h_j}$ holds only if 
$$\epsilon > \omega_{h_j}^2\frac{\Delta t_{\min}^3}{4}$$ 
This case can be ignored in practice since we are interested in small $\epsilon$ which is physically insignificant. Hence, from \eqref{pointwise-bound-case-1}, we have 
\begin{align}
    \hat{h}_{j} \le \delta_{h_j}(\epsilon) = (4\epsilon \omega_{h_j})^{\textstyle\frac{1}{3}}\quad\text{if}\quad\epsilon \le \omega_{h_j}^2\frac{\Delta t_{\min}^3}{4}\label{hj-bnd}  
\end{align}
for $j=1,\ldots,n_h$. We can similarly bound the maximum value of $g_i$, denoted by $\hat{g}_i$, by approximating the area under $\normplus{g_i(t,x(t),u(t))}^2$, for $t\in[t_k,t_{k+1}]$, to obtain 
\begin{align}
    \hat{g}_{i} \le \delta_{g_i}(\epsilon) = (4\epsilon \omega_{g_i})^{\textstyle\frac{1}{3}}\quad\text{if}\quad\epsilon \le \omega_{g_i}^2\frac{\Delta t_{\min}^3}{4}\label{gi-bnd}
\end{align}
for $i=1,\ldots,n_g$. 
%
\section{Control Input Parameterization}\label{ctrl-param-eg}

Let $U = (u_1,\ldots,u_{N_u})$ with $u_k\in\bRnu$, for $k=1,\ldots,N_u$. Then, according to \eqref{aug-ctrl-param}, the control input is given by
\begin{align}
    u(\tau) = \sum_{k=1}^{N_u} \Gamma_u^k(\tau)u_k\label{ctrl-param-expand}
\end{align}
for $\tau\in[0,1]$.
%
\subsection{First-Order-Hold}\label{foh}
Choose $N_u = N$ and  define 
{\small\everymath{\displaystyle}%
\begin{align*}
\Gamma_u^1(\tau) ={} & \left\{ \begin{array}{ll} \frac{\tau_2-\tau}{\tau_2-\tau_1}&\,~\quad\text{if }\tau\in[\tau_1,\tau_2]\\[0.3cm]
0&\,~\quad\text{otherwise}\end{array} \right.\\
\Gamma_u^k(\tau) ={} & \left\{ \begin{array}{ll} \frac{\tau-\tau_{k-1}}{\tau_{k}-\tau_{k-1}} &\,~\text{if }\tau\in[\tau_{k-1},\tau_k]\\[0.5cm]
 \frac{\tau_{k+1}-\tau}{\tau_{k+1}-\tau_{k}} &\,~\text{if }\tau\in[\tau_{k},\tau_{k+1}]\\[0.5cm]
0&\,~\text{otherwise}
\end{array} \right.~~k=2,\ldots,N-1\\
\Gamma_u^N(\tau) ={} & \left\{ \begin{array}{ll} \frac{\tau-\tau_{N-1}}{\tau_{N}-\tau_{N-1}}&\text{if }\tau\in[\tau_{N-1},\tau_N]\\[0.3cm]
0&\text{otherwise}\end{array} \right.
\end{align*}
}%
In other words,
\begin{align}
& u(\tau) = \left(\frac{\tau_{k+1}-\tau}{\tau_{k+1}-\tau_k}\right)u_k + \left(\frac{\tau-\tau_k}{\tau_{k+1}-\tau_k}\right)u_{k+1}
\end{align}
whenever $\tau\in[\tau_k,\tau_{k+1}]$ for some $k=1,\ldots,N-1$. Note that $\Gamma_u^k(\tau_j) = \delta_{jk}$ is the Kronecker delta, for $k,j=1,\ldots,N$. As a result, $u(\tau_k) = u_k$, for $k=1\ldots,N$. 
%
\subsection{Pseudospectral Methods}\label{pseudospectral}
In pseudospectral methods, the control input is parameterized with a basis of orthogonal polynomials, i.e.,
\begin{align}
    \Gamma_u^k(\tau) = \prod_{\genfrac{}{}{0pt}{1}{j=1}{j\ne k}}^{N_u}\frac{2\tau-1 - \eta_j}{\eta_k-\eta_j}\label{lag-interp}
\end{align}
is the Lagrange interpolating polynomial \cite[Sec. 2.5]{davis1975interpolation} of degree $N_u-1$, where $\eta_k\in[-1,1]$, for $k=1,\ldots,N_u$, are related to the roots of orthogonal polynomials (such as members of the Jacobi family \cite[Sec. 10.3]{davis1975interpolation}). For e.g., choosing $\eta_k$ to be the roots of the polynomial
\begin{align}
    \eta \mapsto (1-\eta^2)\dbyd{T_{N_u-1}(\eta)}{\eta}
\end{align}
would correspond to Chebyshev-Gauss-Lobatto (CGL) collocation \cite[Sec. D.1]{malyuta2019discretization}, where $T_{M}:[-1,1] \to \bR$ defined by $\eta \mapsto \cos(M\arccos (\eta) )$ is the Chebyshev polynomial of degree $M$ \cite[Def.  3.3.1]{davis1975interpolation}. The choice of basis functions \eqref{lag-interp} ensures that $\Gamma_u^k\big(\frac{\eta_j+1}{2}\big) = \delta_{jk}$, for each $k,j=1,\ldots,N_u$. As a result, $u\big(\frac{\eta_k+1}{2}\big) = u_k$, for $k=1,\ldots,N_u$.
%
\section{Dynamic Obstacle Avoidance}\label{app:dyn-obstacle-avoid}
We consider a two-dimensional free-final-time path planning problem. The dynamical system describing the vehicle comprises of position $r\in\bR^2$, velocity $v\in\bR^2$, cumulative cost $\check{p} \in \bR$, and acceleration input $u\in\bR^2$. The state $x\in\bR^{5}$, augmented state $\tilde{x}\in\bR^{7}$, and augmented control input $\tilde{u}\in\bR^{3}$ are given by
\begin{subequations}
\begin{align}
x ={} & (r,v,\check{p})\\
\tilde{x} ={} & (x,y,t)\\
\tilde{u} ={} & (u,s)
\end{align}
\end{subequations}
We choose a first-order-hold (see Appendix \ref{foh}) parameterization for the augmented control input. 

The augmented system is given by
\begin{align}
\derv{\tilde{x}} ={} & \left[\begin{array}{c}
     v  \\
     u - c_{\mathrm{d}}\|v\|v \\ 
     \|u\|^2 \\
     \ones{}^\top\normplus{g(t,x,u)}^2 \\
     1
\end{array}\right]\!s
\end{align}
where $c_{\mathrm{d}}\in\bR_+$ is the drag coefficient. The vehicle is subject to the following inequality path constraints: avoidance of $\nobs=10$ moving elliptical obstacles, speed upper bound, and acceleration upper- and lower bounds, which are encoded in the path constraint function $g$ as follows
\begin{align}
g(t, x,u) ={} & \left[ \begin{array}{c} 
1 - \|\check{H}_1(r-\check{q}_1(t))\|^2 \\
\vdots \\
1 - \|\check{H}_{\nobs}(r-\check{q}_{\nobs}(t))\|^2 \\
\|v\|^2 - v^2_{\max} \\ 
\|u\|^2 - u^2_{\max} \\ 
u_{\min}^2 - \|u\|^2 \\ 
\end{array} \right]
\end{align}
The shape matrix and center of the $i$th obstacle at time $t$ are $\check{H}_i$ and $\check{q}_i(t)$, respectively, for $i=1,\ldots,\nobs$. A sinusoidal motion is prescribed for the centers of the obstacles, i.e.,
\begin{align}
    \check{q}_i(t) = \left( \psi_i + \delta\psi_i \sin \left(\theta_i t + \angle\psi_i \right),0 \right)\label{obs-center}
\end{align}
with amplitude $\delta\psi_i$, phase angle $\angle\psi_i$, and frequency $\theta_i$, about the nominal center $\psi_i$. The speed upper bound and the acceleration upper- and lower bounds are $v_{\max}$, $u_{\max}$, and $u_{\min}$, respectively. %
The boundary conditions are specified through function $Q$ as follows
\begin{align}
    Q(\tinit,x(\tinit),\tfinal,x(\tfinal))  ={} & \begin{bmatrix}
        r(\tinit) - \rinit\\
        v(\tinit) - v_{\mathrm{i}}\\
        \check{p}(\tinit)\\
        r(\tfinal) - \rfinal\\
        v(\tfinal) - v_{\mathrm{f}} 
    \end{bmatrix}
\end{align}
where $\rinit$ and $v_{\mathrm{i}}$ are the initial position and velocity, respectively, and $\rfinal$ and $v_{\mathrm{f}}$ are the final position and velocity, respectively. We set the terminal state cost function to $L(\tfinal,x(\tfinal)) = \check{p}(\tfinal)$. Constraint functions $h$ and $P$ are not utilized in this example.

Table \ref{tab:dyn-obstacle-avoidance-param} shows the parameter values chosen for the system and the algorithm to generate the results in Section \ref{subsec:dyn-obstacle-avoid}.  In addition, we choose the following parameter values for the sinusoidal motion given by \eqref{obs-center}
\begin{align*}
    \psi ={} &\!\! 
    \left[ \begin{array}{cccccccccc}
       34 & -32 & 42 & -24 & 34 & -32 &  42 & -24 &  34 & -32 \\
       20 &  20 & 10 &  10 &  0 &   0 & -10 & -10 & -20 & -20
    \end{array} \right]\\    
    \delta\psi ={} & (10)\,\ones{1\times \nobs}\\
    \angle\psi ={} & \frac{\pi}{2} 
    \begin{bmatrix}
       1 & 1 & 0 & 0 & 1 & 1 & 0 & 0 & 1 & 1
    \end{bmatrix}\\
    \theta ={} & \left(\frac{\pi}{20}\right)\ones{1\times\nobs}
\end{align*}
where $\square_i$ is the $i$th column of $\square$, for $\square = \psi,\delta\psi,\angle\psi,\theta$. The scenario with static obstacles in Section \ref{subsec:dyn-obstacle-avoid} is generated by setting $\delta\psi = \zeros{1\times\nobs}$.
\begin{table}[!htpb]
\caption{}\label{tab:dyn-obstacle-avoidance-param}
{\renewcommand{\arraystretch}{1.1}
\begin{tabular}{l|l}
\hline
Parameter & Value\\\hline\\[-0.3cm]
$c_\mathrm{d}$ & $0.01$ m$^{-1}$\\
$\rinit$, $\rfinal$  & $(0,-28)$, $(0,28)$ m\\
$v_{\mathrm{i}}$, $v_{\mathrm{f}}$  & $(0.1,0)$, $(0.1,0)$ m s$^{-1}$\\
$v_{\max}$ & $6$ m s$^{-1}$ \\
$u_{\min}$, $u_{\max}$ & $0.5$, $6$ m s$^{-2}$ \\
$H_i$, $i=1,\dots,\nobs$ & 
$\left[ \begin{array}{cc}
0     & \hphantom{-}0.45\\
0.03  &  0
\end{array} \right]$ \\[0.3cm]
$\tinit$ & $0$ s\\
$\epsilon$ & $10^{-5}$\\
$\mc{U},\mc{S}$ & $\{u\in\bR^{2}\,|\,\|u\|_\infty \le u_{\max}\}$, $[1,60]$\\
$N$ & $10$\\
$\gamma$ & $6.67 \times 10^{3}$\\
$\rho$ & $1.5\times 10^{-3}$\\
\hline
\end{tabular}}
\end{table}
%
\section{6-DoF Rocket Landing}\label{app:6DoF-rocket-landing}
We consider a free-final-time 6-DoF lunar landing scenario based on the formulation given in \cite{szmuk2020successive}, and adopt the notation therein. The state $x\in\bR^{14}$, control input $u\in\bR^3$, augmented state $\tilde{x}\in\bR^{15}$, and augmented control input $\tilde{u}\in\bR^{4}$ are given by
\begin{subequations}
\begin{align}
x ={} & (m,\bm{r}_{\mc{I}},\bm{v}_{\mc{I}},\bm{q}_{\mc{B}\leftarrow\mc{I}},\bm{\omega}_{\mc{B}})\\
u ={} & \bm{T}_{\mc{B}}\\
\tilde{x} ={} & (x,y)\\
\tilde{u} ={} & (\bm{T}_{\mc{B}},s)
\end{align}
\end{subequations}
We choose a first-order-hold (see Appendix \ref{foh}) parameterization for the augmented control input. 

The augmented dynamical system is given by
\begin{align}
\derv{\tilde{x}} = \begin{bmatrix} -\check{\alpha}\|\bm{T}_{\mc{B}}\|\\
    \bm{v}_{\mc{I}}\\
    \frac{1}{m}C_{\mc{I}\leftarrow\mc{B}}(\bm{q}_{\mc{B}\leftarrow\mc{I}})\bm{T}_{\mc{B}} + \bm{g}_{\mc{I}}\\
    \frac{1}{2}\Omega(\bm{\omega}_{\mc{B}})\bm{q}_{\mc{B}\leftarrow\mc{I}}\\
    J_{\mc{B}}^{-1}\left( \bm{r}_{T,\mc{B}}\times\bm{T}_{\mc{B}} - \bm{\omega}_{\mc{B}} \times J_{\mc{B}}\bm{\omega}_{\mc{B}} \right)\\
    \ones{}^\top\normplus{g(t,x,u)}
\end{bmatrix}\!s
\end{align}
The path constraint function $g$ is given by
\begin{align}
    g(t,x,u) ={} & \begin{bmatrix} m_{\mathrm{dry}} - m \\ 
        \| \cot\gamma_{\mathrm{gs}} H_\gamma \bm{r}_{\mc{I}} \|^2 - (\bm{e}_1^\top\bm{r}_{\mc{I}})^2\\ 
        -\bm{e}_1^\top\bm{r}_{\mc{I}}\\
        \|\bm{v}_{\mc{I}}\|^2 - v_{\max}^2\\
        \|2H_\theta\bm{q}_{\mc{B}\leftarrow\mc{I}}\|^2 - (\cos\theta_{\max}-1)^2\\
        \|\bm{\omega}_{\mc{B}}\|^2 - \omega_{\max}^2\\
        \|\cos\delta_{\max}\bm{T}_{\mc{B}}\|^2 - (\bm{e}_1^\top\bm{T}_{\mc{B}})^2\\
        -\bm{e}_1^\top\bm{T}_{\mc{B}}\\
        \|\bm{T}_{\mc{B}}\|^2 - T_{\max}^2\\
        T_{\min}^2 - \|\bm{T}_{\mc{B}}\|^2
     \end{bmatrix}
\end{align}
where the second-order-cone constraints on $\bm{r}_{\mc{I}}$ and $\bm{T}_{\mc{B}}$ are equivalently reformulated to quadratic forms to ensure continuous differentiability \cite[Sec. 3.2.4]{mosekcookbook}. The boundary conditions are specified through function $Q$ as follows
\begin{align}
 Q(\tinit,x(\tinit),\tfinal,x(\tfinal)) ={} & \begin{bmatrix}
    m(\tinit) - m_{\mathrm{wet}}\\
    \bm{r}_{\mc{I}}(\tinit) - \bm{r}_{\mathrm{i}}\\
    \bm{v}_{\mc{I}}(\tinit) - \bm{v}_{\mathrm{i}}\\
    \bm{\omega}_{\mc{B}}(\tinit)\\
    \bm{r}_{\mc{I}}(\tfinal) - \bm{r}_{\mathrm{f}}\\
    \bm{v}_{\mc{I}}(\tfinal) - \bm{v}_{\mathrm{f}}\\
    \bm{q}_{\mc{B}\leftarrow\mc{I}}(\tfinal) - \bm{q}_{\mathrm{id}}\\
    \bm{\omega}_{\mc{B}}(\tfinal)\\    
 \end{bmatrix}
\end{align}
where $\bm{r}_{\mathrm{i}}$ and $\bm{v}_{\mathrm{i}}$ are the initial position and velocity, respectively, $\bm{r}_{\mathrm{f}}$ and $\bm{v}_{\mathrm{f}}$ are the final position and velocity, respectively, and $\bm{q}_{\mathrm{id}}$ is the unit quaternion. The terminal state cost function is set to $L(\tfinal,x(\tfinal)) = -m(\tfinal)$. Constraint functions $h$ and $P$ are not utilized in this example.

Table \ref{tab:rocket-landing-6DoF-param} shows the parameter values chosen for the system and the algorithm to generate the results in Section \ref{subsec:rocket-land-6DoF}. The definitions of $H_{\gamma}$ and $H_{\theta}$ are taken from \cite{szmuk2020successive}.
\begin{table}[!htpb]
\caption{}\label{tab:rocket-landing-6DoF-param}
{\renewcommand{\arraystretch}{1.1}
\begin{tabular}{l|l}
\hline
Parameter & Value \\\hline
$\check{\alpha}$ & $4.53\times 10^{-4}$ s m$^{-1}$\phantom{$X^{X^{X}}$}\\
$\bm{g}_{\mc{I}}$ & $(-1.61,0,0)$ m s$^{-2}$\\
$\bm{r}_{T,\mc{B}}$ & $(-0.25,0,0)$ m\\
$J_{\mc{B}}$ & $\mathrm{diag}(19150,13600,13600)$ kg m$^2$\\
$m_{\mathrm{dry}}$, $m_{\mathrm{wet}}$ & $2100$, $3250$ kg\\
$\bm{r}_{\mathrm{i}}$, $\bm{r}_{\mathrm{f}}$ & $(433,0,250)$, $(10,0,-30)$ m\\
$\bm{v}_{\mathrm{i}}$, $\bm{v}_{\mathrm{f}}$ & $(10,0,-30)$, $(-1,0,0)$ m s$^{-1}$\\
$\gamma_{\mathrm{gs}}$ & $85^\circ$\\
$v_{\max}$ & $50$ m s$^{-1}$\\
$\theta_{\max}$ & $60^\circ$\\
$\omega_{\max}$ & $10^\circ$ s$^{-1}$\\
$\delta_{\max}$ & $45^\circ$ \\
$\tinit$ & $0$ s\\
$T_{\min}$, $T_{\max}$ & $5000$, $22000$ N\\
$\epsilon$ & $10^{-4}$\\
$\mc{U}$, $\mc{S}$ & $\{u\in\bR^{3}\,|\,\|u\|_\infty \le T_{\max}\}$, $[1,60]$\\
$N$ & $5$\\
$\gamma$ & $2\times 10^3$\\
$\rho$ & $2.5\times 10^{-3}$\\
\hline
\end{tabular}}
\end{table}
%
\section{3-DoF Rocket Landing}\label{app:3DoF-rocket-landing}
We consider a convex, fixed-final-time Mars landing scenario based on the dynamical system and constraints in \cite{acikmese2007convex,elango2022customized}. The state $x\in\bR^{7}$ and control input $u\in\bR^{4}$ are given by 
\begin{align}
    x ={} & (r,v,z)\\
    u ={} & (\check{T},\check{\sigma})
\end{align}
where $m = \exp z$ is the vehicle mass, $T = m\check{T}$ is the thrust vector, and $\check{\sigma}$ is an auxiliary control input (slack variable) introduced as a part of the lossless convexification procedure. We choose a zero-order-hold parameterization for the control input. 

We consider the following linear dynamical system
\begin{align}
    \dot{x} ={} & \begin{bmatrix}
        v\\
        \check{T} + (0,0,-g_{\mathrm{m}})\\
        -\check{\alpha}\check{\sigma}
    \end{bmatrix}\\ 
    ={} & \underbrace{\begin{bmatrix}
       \zeros{3\times 3} & \eye{3} & \zeros{3\times 1} \\ 
         & \zeros{3\times 7} & \\
         & \zeros{1\times 7} &  
    \end{bmatrix}}_{A}x + \underbrace{\begin{bmatrix} \zeros{3\times 4} \\ \begin{array}{cc} \eye{3} & \zeros{3\times 1} \\ \zeros{1\times 3} & -\check{\alpha} \end{array} \end{bmatrix}}_B u + \underbrace{\begin{bmatrix}
        \zeros{5\times 1} \\ -g_{\mathrm{m}}\\ 0 
    \end{bmatrix}}_w\nonumber
\end{align}
The path constraint function $g$ is given by 
\begin{align}
    & g(t,x,u) ={} \\
    & \begin{bmatrix}
      \|\cot\gamma_{\mathrm{gs}}E_{12}r\|^2 - (e_3^\top r)^2\\
      -(e_3^\top r)\\
      \|v\|^2 - v_{\max}^2\\
      z - z_1(t)\\
      -z + \max\{\log m_{\mathrm{dry}},z_0(t)\}\\
     -e_3^\top\check{T}+\check{\sigma}\cos\theta_{\max}\\
     \|\check{T}\|^2 - \check{\sigma}^2\\
     -\check{\sigma}\\
     \mu_{\min}(t)\left(1-(z-z_0(t))+\frac{1}{2}(z-z_0(t))^2\right) - \check{\sigma}\\
    -\mu_{\max}(t)(1-(z-z_0(t))) + \check{\sigma}
    \end{bmatrix}\nonumber
\end{align}
where $E_{12} = [\eye{2}~\zeros{2\times 1}]$, $e_3 = (0,0,1)$, and 
\begin{subequations}
\begin{align}
    z_0(t) ={} & \log (m_{\mathrm{wet}} - \check{\alpha}T_{\max}t)\\
    z_1(t) ={} & \log (m_{\mathrm{wet}} - \check{\alpha}T_{\min}t)\\
    \mu_{\min}(t) ={} & T_{\min}\exp(-z_0(t))\\
    \mu_{\max}(t) ={} & T_{\max}\exp(-z_0(t))
\end{align}    
\end{subequations}
Lossless convexification holds if $\|T(t)\| = \sigma(t) = m(t)\check{\sigma}(t)$ $\ae$ $t\in\tspan$.
The boundary conditions are specified through function $Q$ as follows
\begin{align}
    Q(\tinit,x(\tinit),\tfinal,x(\tfinal)) ={} \begin{bmatrix}
        r(\tinit) - \rinit\\
        v(\tinit) - v_{\mathrm{i}}\\
        z(\tinit) - \log m_{\mathrm{wet}}\\
        r(\tfinal) - \rfinal\\
        v(\tfinal) - v_{\mathrm{f}}
    \end{bmatrix}
\end{align}
where $r_{\mathrm{i}}$ and $v_{\mathrm{i}}$ are the initial position and velocity, respectively, and $r_{\mathrm{f}}$ and $v_{\mathrm{f}}$ are the final position and velocity, respectively. The terminal state cost function is set to $L(\tfinal,x(\tfinal)) = -z(\tfinal)$. The running cost function $Y$, and constraint functions $h$ and $P$, are not utilized in this example.

Table \ref{tab:rocket-landing-3DoF-param} shows the parameter values chosen for the system and the algorithm to generate the results in Section \ref{subsec:rocket-land-3DoF}.

\begin{table}[!htpb]
\caption{}\label{tab:rocket-landing-3DoF-param}
{\renewcommand{\arraystretch}{1.1}
\begin{tabular}{l|l}
\hline
Parameter &  Value \\ \hline 
$\check{\alpha}$ &   $4.53\times 10^{-4}$ s m$^{-1}$\phantom{$X^{X^{X}}$}\\
$g_{\mathrm{m}}$ & $3.71$ m s$^{-2}$\\
$\rinit$, $\rfinal$ & $(2000,0,1500)$, $(0,0,0)$ m\\
$v_{\mathrm{i}}$, $v_{\mathrm{f}}$ & $(80,30,-75)$, $(0,0,0)$ m s$^{-1}$\\
$m_{\mathrm{dry}}$, $m_{\mathrm{wet}}$ & $1505$, $1905$ kg\\
$\gamma_{\mathrm{gs}}$ & $84^\circ$\\
$v_{\max}$ & 139 m s$^{-1}$\\
$\theta_{\max}$ & $40^\circ$\\
$T_{\min}$, $T_{\max}$ & $4971.6$, $13258$ N\\
$[\tinit,\tfinal]$ & $[0,84]$ s\\
$\epsilon$ & $10^{-5}$\\
$\mc{U}$ & $\{u\in\bR^4\,|\,\|u\|_\infty \le T_{\max}\}$\\
$N$ & $8$ \\
$\gamma$ & $6.67 \times 10^3$ \\
$\rho$ & $1.5\times 10^{-3}$\\
\hline
\end{tabular}}
\end{table}
%
\end{document}